\documentclass[10pt]{article}
\usepackage[utf8]{inputenc}
\usepackage[T1]{fontenc}
\usepackage[english]{babel}
\usepackage{cancel}
\usepackage{bm}
\usepackage[table,dvipsnames]{xcolor}
\definecolor{blun}{cmyk}{0.8, 0.5, 0, 0.7} % blunotte
\usepackage{amsmath} % package qui change et ameliore les commandes pour les formules
\usepackage{amsfonts} % package pour caracteres mathematiques supplementaires
\usepackage{amssymb} % package pour symboles mathematiques supplementaires
\usepackage{amsthm} % verison amelioree de \newtheorem + autre options
\usepackage{mathrsfs} % Police calligraphique
\usepackage{pifont} % symboles sympa a l'aide de '\ding{numero-du-symbole}'
\usepackage{verbatim} % for comment environment
\usepackage[bookmarks,colorlinks,pdfauthor={Cinzia Soresina},linkcolor=blun,citecolor=blun]{hyperref}
\usepackage{enumerate} % for different enumerators
\usepackage{graphicx}
\usepackage[TABBOTCAP]{subfigure}
\DeclareMathAlphabet{\mathscrbf}{OMS}{mdugm}{b}{n}

\usepackage{geometry}
\newcommand{\N}  {\mathbb{N}}
\newcommand{\R}  {\mathbb{R}}

\newcommand{\var} {\varepsilon}
\newcommand{\pa} {\partial}
\newcommand{\tr}  {\textnormal{tr }}

\newcommand{\beq}{\begin{equation}}
\newcommand{\eeq}{\end{equation}}
\newcommand{\ben}{\begin{eqnarray}}
\newcommand{\een}{\end{eqnarray}}
\newcommand{\beno}{\begin{eqnarray*}}
\newcommand{\eeno}{\end{eqnarray*}}
\newcommand{\beqno}{\begin{equation*}}
\newcommand{\eeqno}{\end{equation*}}
\newcommand{\be}{\begin{enumerate}}
\newcommand{\ee}{\end{enumerate}}
\newcommand{\bitem}{\begin{itemize}}
\newcommand{\eitem}{\end{itemize}}
\newcommand{\baln}{\begin{align}}
\newcommand{\ealn}{\end{align}}
\newcommand{\balno}{\begin{align*}}
\newcommand{\ealno}{\end{align*}}

\newtheorem{thm}{Theorem}[section]

\usepackage{amsbsy}

\begin{document}

\begin{center}
{\bf{\large{About reaction-diffusion systems involving the Holling-type II and the Beddington-DeAngelis functional responses for predator-prey models}}}
\end{center}
\bigskip

\bigskip

\begin{center}
F. CONFORTO\footnotemark[1], L. DESVILLETTES\footnotemark[2], AND C. SORESINA\footnotemark[3]
\end{center}

\noindent\footnotetext[1]{\label{adr1}Dipartimento di Scienze Matematiche e Informatiche, Scienze Fisiche e Scienze della Terra, Università di Messina, Viale Stagno d'Alcontres 31, 98166 Messina, Italy}
\noindent\footnotetext[2]{\label{adr2} Universit\'{e} Paris Diderot, Sorbonne Paris Cit\'{e}, Institut de Math\'{e}matiques de Jussieu-Paris Rive Gauche, UMR 7586, CNRS, Sorbonne Universit\'{e}s, UPMC Univ. Paris 06, F-75013, Paris, France}
\noindent\footnotetext[3]{\label{adr3} Centro de Matem\'atica, Aplica\c{c}\~{o}es Fundamentais e Investiga\c{c}\~ao Operacional, Universidade de Lisboa, Faculty of Science, Campo Grande, 1749-016 Lisboa, Portugal}

% \section{Introduction}\label{chap1}

\begin{abstract}
We consider in this paper a microscopic model (that is, a system of three reaction-diffusion equations) incorporating the dynamics of handling and searching predators, and show that its solutions converge when a small parameter tends to $0$ towards the solutions of a reaction-cross diffusion system of predator-prey type involving a Holling-type II or Beddington-DeAngelis functional response. We also provide a study of the Turing instability domain of the obtained equations and (in the case of the Beddington-DeAngelis functional response) compare it to the same instability domain when the cross diffusion is replaced by a standard diffusion.
\end{abstract}
\bigskip

\bigskip 

\bigskip

\bigskip 

\bigskip 

\section{Introduction}\label{s1}
\subsection{General presentation}
Complex functional responses are widely used in predator-prey models \cite{Holl, Ivlev, bedd, DeAng2, abragi}. For example, the Holling-type II functional response \cite{Holl} is based on the idea that predators will catch a limited proportion of available prey in the case when preys are abundant. Denoting with $N:=N(t)$ the prey biomass, and with $P:=P(t)$ the predator biomass, this type of functional response leads to the following set of two ODEs:
\begin{align}\begin{split}
\dot{N}&=r_0g(N)-\dfrac{b NP}{1+k N},\\
\dot{P}&=\dfrac{cb NP}{1+k N}-\mu P,
\end{split}\label{ODEppHoll}\end{align}
where $r_0,\;b,\;c,\;k,\;\mu>0$, and where the function $g$ describes the prey growth and can be either linear, that is $g(N)=N$, or involve a logistic part as $g(N)=(1-\eta N)N$ with $\eta > 0$ \cite{Bazykin}. Note that when $g(N)=N$ and $k=0$, one recovers the classical  Lotka-Volterra predator-prey model.

If one also wishes to take into account the competition between predators when they try to catch prey, the slightly more complex Beddington-DeAngelis functional response can be introduced \cite{bedd,DeAng2}:
\begin{align}\begin{split}
\dot{N}&=r_0g(N)-\dfrac{b NP}{1+k N+h P},\\
\dot{P}&=\dfrac{cb NP}{1+k N+h P}-\mu P,
\end{split}\label{ODEppBdA}\end{align}
where $r_0,\;b,\;c,\;k,\;h>0$.

An important point in the sequel will be the observation that predator-prey models with the Beddington-DeAngelis functional response are known to produce patterns (coming out of a Turing instability) when diffusion terms with suitable rates (denotes by $d_N,\;d_P$) are added to the reaction term \cite{zhang2012, Haque2012}. In such a situation, the system becomes: 
\begin{align}\begin{split}
\partial_t N-d_N \Delta_x N&=r_0g(N)-\dfrac{b NP}{1+k N+h P},\\
\partial_t P-d_P \Delta_x P&=\dfrac{cb NP}{1+k N+h P}-\mu P.
\end{split}\label{PDEppBdAlin}\end{align}
On the other hand, no patterns are known to appear in the case of a reaction-diffusion predator-prey model with a Holling-type II functional response and diffusion terms as in (\ref{PDEppBdAlin}) \cite{Alonso2002} (patterns may however appear when richer dynamics are considered, for example when one adds quadratic intra-predator interaction or fighting term \cite{mcgehee2005,malchow2008}, or a density-dependent predator mortality \cite{malchow2008}). In all these cases, a fundamental assumption is that the diffusion coefficients of the two species must be different;  patterns appearing taking equal diffusion coefficient are studied in \cite{petrovskii1999,petrovskii2001,medvinsky2002}).

Several works were written in the past in order to obtain a derivation of 
the Holling-type II and Beddington-DeAngelis functional responses out of
simple and realistic 
%An interesting modeling issue consists in finding a realistic
 ``microscopic models'' 
 %(in terms of time scales)
  which in some limit lead, at least formally, to \eqref{ODEppHoll} or \eqref{ODEppBdA}. 
 Such a model was designed by Metz and Diekmann \cite{metz2014} for the Holling-type II functional response, and by Geritz and Gyllenberg \cite{geritz2012}, Huisman and De Boer \cite{Huisman1997}, for the Beddington-DeAngelis one, and references therein.
 
Metz and Diekmann proposed a system of three ODEs, in which the predators are divided in two classes (respectively called searching and handling), while the interaction between predators and preys is treated in a quite simple way (standard Lotka-Volterra terms are used). Predators which are searching for preys become handling with a rate proportional to the number of preys and come back to the searching state with a constant rate. Only handling predators contribute to the reproduction (and give birth to a searching predator), while the mortality rate (in absence of prey) is constant and equal for the two classes. The searching-handling switch is supposed to happen on a much faster time scale than the reproduction and mortality processes. The corresponding parameter in the system of ODEs is therefore called $1/\var$, and the system writes (for some $r_0$, $\alpha$, $\tilde{\gamma}$, $\Gamma$, $\mu$, $\var >0$)
\begin{align}\begin{split}
\dot{N}^\varepsilon&=r_0 N^\varepsilon-\alpha N^\varepsilon p_s^\varepsilon,\\
\dot{p}_s^\varepsilon&=\dfrac{1}{\varepsilon}\left(-\alpha N^\varepsilon p_s^\varepsilon + \tilde{\gamma} p_h^\varepsilon \right) +\Gamma p_h^\varepsilon - \mu p_s^\varepsilon,\\
\dot{p}_h^\varepsilon&=\dfrac{1}{\varepsilon}\left(\alpha N^\varepsilon p_s^\varepsilon - \tilde{\gamma} p_h^\varepsilon\right)- \mu p_h^\varepsilon,
\label{sysOdoODE3}
\end{split}\end{align}
where $N^\var$ still is the density of preys, while $p_h^\varepsilon$ and $p_s^\varepsilon$ are the respective densities of handling and searching predators.
It is shown that in the formal limit $\varepsilon\to 0$, one gets $N^\varepsilon\to N$ and $p_s^\varepsilon+p_h^\varepsilon\to P$ where $N,\;P$ satisfy \eqref{ODEppHoll} (with $g(N):= N$, $b :=\alpha,\; k :=\alpha/\tilde{\gamma},\;c :=\Gamma/\tilde{\gamma}$) \cite{metz2014}.

A similar procedure was later applied by Geritz and Gyllenberg in a more complex situation \cite{geritz2012}: they divided not only the predator population into searchers and handlers, but also structured the prey population into two classes, the class of active preys (typically foraging) and prone to predation, and the class of those prey individuals who have found a refuge and cannot be caught by predators. In this way, they derived the Beddington-DeAngelis functional response in terms of mechanisms at the individual level avoiding the usual interference between predators. Previously, Huisman and De Boer \cite{Huisman1997}, starting from a different four-dimensional model 
%(prey, predators, an interaction complex of prey and predator, handling predators),
% use a reaction schemes and quasi-steady-state assumptions,
also obtained a system of two ordinary differential equations (they however simplified a complicated quadratic expression with a Pad\'e approximation to recover the standard formula of the Beddington-DeAngelis functional response). In both cases, two different time scales were exploited.
\medskip

We are interested in this paper in the introduction of diffusion processes in the asymptotic problem \eqref{sysOdoODE3}. Denoting by $d_1,\;d_2,\;d_3>0$ the diffusion rates of preys, searching predators and handling predators respectively, one can write, keeping the reaction term of (\ref{sysOdoODE3}), the following reaction-diffusion system:
\begin{align}\begin{split}
\pa_t N^\var-d_1 \Delta_x N^\var&=r_0\left(1-\eta N^\var\right) N^\var-\alpha N^\varepsilon p_s^\varepsilon,\\  
\pa_t p_s^\var-d_2\Delta_x p_s^\var&=\dfrac{1}{\varepsilon}\left(-\alpha N^\varepsilon p_s^\varepsilon + \tilde{\gamma}\, p_h^\varepsilon \right) +\Gamma p_h^\varepsilon - \mu p_s^\varepsilon,\\
\pa_t p_h^\var - d_3 \Delta_x p_h^\var &= \dfrac{1}{\varepsilon}\left(\alpha N^\varepsilon p_s^\varepsilon - \tilde{\gamma}\, p_h^\varepsilon\right)- \mu p_h^\varepsilon.
\end{split}\label{PDEHoll}\end{align}
Note that we systematically expect the diffusion rate $d_3$ of handling predators to be smaller than the diffusion rate of searching predators $d_2$. The formal limit of this system when $\var \to 0$ is the set of two reaction cross-diffusion equations: 
\begin{align}\begin{split}
\pa_t& N - d_1 \Delta_x N = r_0(1-\eta N) N - \frac{ \tilde{\gamma}\alpha N}{\alpha N + \tilde{\gamma}}  P,\\
\pa_t& P -  \Delta_x \bigg( \frac{ d_2 \tilde{\gamma} + d_3 \alpha N}{\alpha N + \tilde{\gamma}} \,P \bigg) = 
 \Gamma\frac{\alpha N\, P}{\alpha N + \tilde{\gamma}} -  \mu P  ,
\end{split}\label{limHollcross}\end{align}
in which the reaction terms are identical (up to the change of name of the constant parameters) to those of \eqref{ODEppHoll}, but in which the diffusion term relative to predators is much more complicated than a constant times Laplacian of $P$ (terms like $d\Delta_x P$ will be systematically called linear diffusive terms in the sequel, while cross diffusion refers to terms like $\Delta_x (f(N,P)\,P)$, where $f$ is a smooth non-constant function of $N,\: P$, as in the second equation of \eqref{limHollcross}). It can be noticed that the resulting cross-diffusion term is a convex combination of the diffusion coefficients $d_2$ and $d_3$ of the microscopic system. In Subsection \ref{RigRef} of the Introduction of this paper, we state a rigorous theorem showing that convergence of solutions to system \eqref{PDEHoll} towards solutions to system \eqref{limHollcross} indeed holds when suitable functional spaces are introduced. 

The same procedure can be applied in the case of Beddington-DeAngelis like functional response, that is a system of ODEs close to \eqref{ODEppBdA}. First, we introduce a system of three ordinary differential equations modeling the interaction between preys, handling and searching predators as in \eqref{PDEHoll}, but in which we also take into account the competition among predators when they look for preys. This is done thanks to the introduction of the denominator $1+\xi p_s$, for some $\xi>0$,  in the interaction term between predators and preys. The system writes as follows: 
\begin{align}\begin{split}
\dot{N}^\varepsilon&=r_0 (1-\eta N^\varepsilon)N^\var-\dfrac{\alpha N^\varepsilon p_s^\varepsilon}{1+\xi p_s^\varepsilon},\\
\dot{p}_s^\varepsilon&=\dfrac{1}{\varepsilon}\left(-\dfrac{\alpha N^\varepsilon p_s^\varepsilon}{1+\xi p_s^\varepsilon} + \tilde{\gamma} p_h^\varepsilon \right) +\Gamma p_h^\varepsilon - \mu p_s^\varepsilon,\\
\dot{p}_h^\varepsilon&=\dfrac{1}{\varepsilon}\left(\dfrac{\alpha N^\varepsilon p_s^\varepsilon}{1+\xi p_s^\varepsilon} - \tilde{\gamma} p_h^\varepsilon\right)- \mu p_h^\varepsilon.
\end{split}\label{sysBdAODEeps}\end{align}
Its formal limit when $\var\to0$ is then a system close to \eqref{ODEppBdA}, also obtained in \cite{Huisman1997} starting from a
% more complex 
system of four ODEs in which all interactions are linear/quadratic.

A reaction-diffusion system corresponding to \eqref{sysBdAODEeps}, where
the diffusion of preys, searching predators, and handling predators is taken
into account through diffusion rates $d_1,\;d_2,\;d_3>0$, writes: 
\begin{align}\begin{split}
\pa_t N^\var-d_1 \Delta_x N^\var&=r_0\left(1-\eta N^\var\right) N^\var-\dfrac{\alpha N^\var p_s^\var}{1 + \xi p_s^\var}, \\  
\pa_t p_s^\var-d_2\Delta_x p_s^\var&=\dfrac{1}{\var}\bigg(- \frac{\alpha N^\var p_s^\var}{1 + \xi p_s^\var} + \tilde{\gamma} p_h^\var \bigg)  + \Gamma p_h^\var - \mu p_s^\var,\\
\pa_t p_h^\var - d_3 \Delta_x p_h^\var &= -\frac1{\var} \bigg( - \frac{\alpha N^\var p_s^\var}{1 + \xi p_s^\var} + \tilde{\gamma} p_h^\var \bigg)  - \mu p_h^\var. 
\end{split}\label{sysBdAPDEeps}\end{align}
We present in Subsection \ref{RigRef} of the Introduction a rigorous result of convergence of the solutions to this system towards the solution of a reaction-cross diffusion system where the reaction part is close to \eqref{ODEppBdA}. This system writes
%A=\gamma+N-P,\quad B=\gamma+N+P,\quad \Delta=(\gamma+N+P)^2-4NP.\label{ABDelta}
\begin{align}\begin{split}
\dfrac{\partial N}{\partial T}& - d_1\Delta_x N = r\left(1-\eta\,N\right)N-\dfrac{ 2\alpha \tilde{\gamma} NP}{\tilde{\gamma}+ \alpha N+ \tilde{\gamma}\xi P+\sqrt{(\tilde{\gamma}+ \alpha N - \tilde{\gamma}\xi P)^2 + 4 \tilde{\gamma}2\,\xi P}}, \\
\dfrac{\partial P}{\partial T}& - \Delta_x \left(f(N,P)P\right)=\dfrac{2\alpha \Gamma NP}{\tilde{\gamma}+ \alpha N+ \tilde{\gamma}\xi P+\sqrt{(\tilde{\gamma}+ \alpha N - \tilde{\gamma}\xi P)^2 + 4 \tilde{\gamma}2\,\xi P}}-\mu P,
\end{split}\label{sysBdAPDElim}\end{align}
where 
%\begin{multline*}
\begin{equation}\label{suppl}
f(N,P)=d_2\dfrac{2 \tilde{\gamma} }{\tilde{\gamma}+ \alpha N - \tilde{\gamma}\xi P+\sqrt{(\tilde{\gamma}+ \alpha N - \tilde{\gamma}\xi P)^2 + 4 \tilde{\gamma}2\,\xi P}}\\+d_3\dfrac{2 \alpha N}{\tilde{\gamma}+ \alpha N+ \tilde{\gamma}\xi P+\sqrt{(\tilde{\gamma}+ \alpha N - \tilde{\gamma}\xi P)^2 + 4 \tilde{\gamma}2\,\xi P}}.
\end{equation}
%\end{multline*}

The proof of this result as well as the proof of the theorem corresponding to the Holling-type II functional response, is based on estimates coming out of two classes of methods. On one hand, we use the duality lemmas devised for reaction-diffusion systems by M. Pierre and D. Schmitt \cite{DS}. More precisely, we use an improved version of those lemmas allowing to recover
 $L^p$ bounds for $p>2$ for the solutions of such systems \cite{CDF,BDF}. On the other hand, we also use entropy-like functionals which are strongly reminiscent of those used in works in which microscopic models for the Shigesada-Kawasaki-Teramoto system \cite{SKT} are studied \cite{DeTr}.
 
These proofs can also be compared to recent results in which reaction-cross diffusion systems are obtained as limits of standard reaction-diffusion systems with more equations, in the context of chemistry or biology \cite{mu, izmi, code, BP1, BP2, HMN, HHP1, HHP2, Rol}.  
\bigskip

As already mentioned, the Beddington-DeAngelis like functional response is particularly interesting since it is known that predator-dependent functional responses can lead to patterns when (linear) diffusion terms are added to the reaction terms, like for example in the following system (with reaction terms identical to those in (\ref{sysBdAPDElim})):

\begin{align}\begin{split}
\dfrac{\partial N}{\partial T}& - D_1\Delta_x N = r\left(1-\eta\,N\right)N-\dfrac{ 2\alpha \tilde{\gamma} NP}{\tilde{\gamma}+ \alpha N+ \tilde{\gamma}\xi P+\sqrt{(\tilde{\gamma}+ \alpha N - \tilde{\gamma}\xi P)^2 + 4 \tilde{\gamma}2\,\xi P}}, \\
\dfrac{\partial P}{\partial T}& - D\Delta_x P =\dfrac{2\alpha \Gamma NP}{\tilde{\gamma}+ \alpha N+ \tilde{\gamma}\xi P+\sqrt{(\tilde{\gamma}+ \alpha N - \tilde{\gamma}\xi P)^2 + 4 \tilde{\gamma}2\,\xi P}}-\mu P,
\end{split}\label{sysBdAPDElin}\end{align}

However, if one consider that the Beddington-DeAngelis like functional response is coming out of an asymptotics when $\var \to 0$ of \eqref{sysBdAPDEeps}, one should rather study the possible appearance of patterns starting from system \eqref{sysBdAPDElim}. We will describe the study that we performed concerning this issue in Subsection \ref{TIa intro} of the Introduction. Note that for Holling-type II functional response, no patterns appear when the cross diffusion model \eqref{limHollcross} is considered, at least under the (biologically reasonable) assumption $d_3<d_2$, so that the qualitative behavior of (\ref{PDEppBdAlin}) and (\ref{limHollcross}) is not expected to be different.

%%%%%%%%%%%%%%%%%%%%%%%%%%%%%%%%%%%%%%%%%%%%%%%%%%%%%%%%%%%%%%%%%%%%%%%%%%%%%%%%%%%%%%%%%%%%%%%%%%%%%%%%%%%%%%%%%%%%%%%%%%%%%%%%%%%%%%%%%
%Scaling!!!
%\begin{equation}\begin{cases}
%\dfrac{\partial n}{\partial t} - d_1 \Delta_x n     = r_0\left(1-\dfrac{n}{n_{max}}\right)n- \dfrac{\alpha n p_s}{1+\xi p_s}\\[0.3cm]
%\dfrac{\partial p_s}{\partial t} - d_2 \Delta_x p_s=\left(-\dfrac{\alpha n p_s}{1+\xi p_s}+\gamma p_h\right)+\gamma \epsilon_0 p_h-\mu p_s\\[0.3cm]
%\dfrac{\partial p_h}{\partial t} - d_3 \Delta_x p_h=\left(\dfrac{\alpha n p_s}{1+\xi p_s}-\gamma p_h\right)-\mu p_h
%\label{fullmodel2}\end{cases}\end{equation}
%Introducing $N\var=\epsilon n,\; k=\epsilon n_{max},\; \tilde{\gamma}=\epsilon \gamma,\; \Gamma=\gamma \epsilon_0$ in the full model \eqref{fullmodel2}, we obtain
%\begin{equation}\begin{cases}
%\dfrac{\partial N}{\partial t} - d_1 \Delta_x N     = r_0\left(1-\dfrac{N}{k}\right)N- \dfrac{\alpha N p_s}{1+\xi p_s}\\[0.3cm]
%\dfrac{\partial p_s}{\partial t} - d_2 \Delta_x p_s=\dfrac{1}{\epsilon}\left(-\dfrac{\alpha N p_s}{1+\xi p_s}+\tilde{\gamma} p_h\right)+\Gamma p_h-\mu p_s\\[0.3cm]
%\dfrac{\partial p_h}{\partial t} - d_3 \Delta_x p_h=-\dfrac{1}{\epsilon}\left(-\dfrac{\alpha N p_s}{1+\xi p_s}+\tilde{\gamma} p_h\right)-\mu p_h
%\end{cases}\label{qssasystem2}\end{equation}

\subsection{Rigorous results of convergence}\label{RigRef}
We consider in this subsection the system (\ref{sysBdAPDEeps}) in a  smooth bounded open subset 
$\Omega$ of $\R^d$: 
\begin{align}
\pa_t N^\var-d_1 \Delta_x N^\var&=r_0\left(1-\eta N^\var\right) N^\var-\dfrac{\alpha N^\var p_s^\var}{1 + \xi p_s^\var}, \label{eqN1bis}\\  
\pa_t p_s^\var-d_2\Delta_x p_s^\var&=\dfrac{1}{\var}\bigg(- \frac{\alpha N^\var p_s^\var}{1 + \xi p_s^\var} + \tilde{\gamma} p_h^\var \bigg)  + \Gamma p_h^\var - \mu p_s^\var, \label{eqN2bis}\\
\pa_t p_h^\var - d_3 \Delta_x p_h^\var &= -\frac1{\var} \bigg( - \frac{\alpha N^\var p_s^\var}{1 + \xi p_s^\var} + \tilde{\gamma} p_h^\var \bigg)  - \mu p_h^\var, \label{eqN3bis}
\end{align}
together with homogeneous Neumann boundary conditions ($\hat{n}(x)$
denoting the exterior normal to $\Omega$ at a point $x\in\pa\Omega$)
\begin{equation}\label{NBC}
\forall x \in \pa\Omega, \qquad \hat{n}(x)\cdot \nabla_x N^\var=0,\quad \hat{n}(x)\cdot \nabla_x p_s^\var=0, \quad \hat{n}(x)\cdot \nabla_x p_h^\var=0 . 
\end{equation}
 All parameters $r_0$, $\alpha$, $\tilde{\gamma}$, $\Gamma$, $\mu$, $d_i$
 in this system are strictly positive, except $\eta$ and $\xi$, which are 
 supposed to be nonnegative.
 When $\eta=0$, no direct logistic saturation is imposed on the preys, while when $\xi=0$, no competition between predators is assumed. Note 
 that when both $\eta$ and $\xi$ are equal to zero, the reaction part of the system \eqref{eqN1bis}-\eqref{eqN3bis} reduces to \eqref{sysOdoODE3}.
\bigskip

We begin with a rigorous result for the passage to the limit $\var\to 0$ in the case $\xi=0$:
\begin{thm}\label{thmsimple}
Let $\Omega$ be a smooth bounded open subset  of $\R^d$ (for some dimension $d \in \N - \{0\}$),  $d_1,d_2,d_3>0$ be diffusion rates, $r_0,\; \alpha,\; \tilde{\gamma}\;, \Gamma,\; \mu>0$ and $\eta\geq0$ be parameters, and $N_{in}:=N_{in}(x) \ge 0$, $p_{h, in} := p_{h, in}(x) \ge 0$, $p_{s, in} := p_{s, in}(x) \ge 0$ be nonnegative initial data respectively in $L^{\infty}(\Omega)$, $L^{2+\delta}(\Omega)$, and $L^{2+\delta}(\Omega)$ for some $\delta>0$. 

Then for each $\var>0$, there exists a unique global classical (for $t>0$) solution ($N^\var$, $p_h^\var$, $p_s^\var$) to system \eqref{eqN1bis} -- \eqref{NBC} with $\xi=0$
% and satisfying homogeneous Neumann boundary conditions (and 
(with the initial data defined above).

Moreover, when $\var \to 0$, one can extract from $N^\var$ a subsequence which is bounded in $L^{\infty}([0,T] \times\Omega)$ for all $T>0$ and converges a.e. towards a function $N\ge 0$ lying in  $L^{\infty}([0,T] \times\Omega)$. One can also extract from $p_s^\var$ (resp. $p_h^\var$) a subsequence which converges weakly in  $L^{2+\delta}([0,T] \times\Omega)$ towards a function $p_s\ge 0$ (resp. $p_h\ge 0$) lying in $L^{2+\delta}([0,T] \times\Omega)$ for all $T>0$ and some $\delta>0$.
% Finally, one can extract from $p_h^\var$ a subsequence which converges weakly in  $L^{2+\delta}([0,T] \times\Omega)$  towards a function $p_h\ge 0$ lying in  $L^{2+\delta}([0,T] \times\Omega)$  for all $T>0$ and some $\delta>0$.

Finally, $N$, $p_s$ and $p_h$ are very-weak solutions of the cross-diffusion system
\begin{align}
\pa_t& N - d_1 \Delta_x N= r_0 (1-\eta N) N - \alpha N p_s, \label{eqF1}\\
\pa_t& (p_s + p_h)  -  \Delta_x (d_2 p_s + d_3 p_h) =  \Gamma p_h - \mu (p_h + p_s),\label{eqF2} \\[0.1cm]
&\alpha N p_s = \tilde{\gamma} p_h , \label{eqF3}
\end{align}
together with the homogeneous Neumann boundary conditions  
\begin{equation}\label{NBC2}
\forall x \in \pa\Omega, \qquad \hat{n}(x)\cdot \nabla_x N=0,\quad \hat{n}(x)\cdot \nabla_x p_s=0, 
\end{equation}
and with initial data 
\begin{equation}\label{init}
N(0,x) = N_{in}(x),\qquad (p_s + p_h)(0,x) = p_{s,in}(x)+p_{h,in}(x),
\end{equation}
in the following sense: identity (\ref{eqF3}) holds a.e.,
 and for all
 $\varphi,\psi \in \mathcal{C}^2_c([0,+\infty)\times\bar{\Omega})$
 such that $ \nabla_x \varphi \cdot n_{|_{\partial \Omega}}=0,\; \nabla_x \psi \cdot n_{|_{\partial \Omega}}=0$,  
%\begin{multline*}
\begin{equation}\label{vw1}
-\int_{0}^{\infty}\int_{\Omega} N\partial_t\varphi-\int_{\Omega}\varphi(0,\cdot)N_{in}-d_1\int_{0}^{\infty}\int_{\Omega}N\Delta_x\varphi%\\
=\int_{0}^{\infty}\int_{\Omega}(r_0 (1-\eta N) N - \alpha N p_s)\varphi,
\end{equation}
%\end{multline*}
%\begin{multline*}
\begin{equation}\label{vw2}
-\int_{0}^{\infty}\int_{\Omega} (p_s + p_h)\partial_t\psi
-\int_{\Omega}\psi(0,\cdot)(p_{s,in}+p_{h,in})%\\
-\int_{0}^{\infty}\int_{\Omega}(d_2 p_s + d_3 p_h)\Delta_x\psi
=\int_{0}^{\infty}\int_{\Omega}(\Gamma p_h - \mu (p_h + p_s))\psi .
 \end{equation}
%\end{multline*}
%$$\forall\varphi,\psi \in \mathcal{C}^2_C([0,+\infty)\times\bar{\Omega}), \nabla_x \varphi \cdot n_{|_{\partial \Omega}}=0,\; \nabla_x \psi \cdot n_{|_{\partial \Omega}}=0.$$

Note that the reaction-cross diffusion system (\ref{eqF1}) - (\ref{eqF3}) together with the homogeneous boundary condition (\ref{NBC2})   can be rewritten in the simpler form (with $P=p_h + p_s$)
\begin{align}
\pa_t& N - d_1 \Delta_x N = r_0(1-\eta N) N - \frac{ \tilde{\gamma}\alpha N}{\alpha N + \tilde{\gamma}}  P,\label{sfcrossS1}\\
%\pa_t& P -  \Delta_x \bigg( \frac{ d_2 \tilde{\gamma} + d_3 \alpha N}{\alpha N + \tilde{\gamma}} \,P \bigg) = \bigg[ \frac{\alpha N}{\alpha N + \tilde{\gamma}} (\Gamma - \mu) -  \frac{\tilde{\gamma}}{\alpha N + \tilde{\gamma}} \mu \bigg]\, P,
\pa_t& P -  \Delta_x \bigg( \frac{ d_2 \tilde{\gamma} + d_3 \alpha N}{\alpha N + \tilde{\gamma}} \,P \bigg) = 
 \left(\Gamma\frac{\alpha N}{\alpha N + \tilde{\gamma}} -  \mu \right)P,\label{sfcrossS2}
\end{align}
%together with homogeneous Neumann boundary conditions 
\begin{equation}\label{newNBC}
\forall x \in \pa\Omega, \qquad \hat{n}(x)\cdot \nabla_x N=0,\quad \hat{n}(x)\cdot \nabla_x P=0,
\end{equation}
with initial data 
\begin{equation}\label{newint}
N(0,x) = N_{in}(x),\qquad P(0,x) = P_{in}(x) := p_{s,in}(x)+p_{h,in}(x).
\end{equation}
%and its very-weak formulation is
%\begin{align*}
%&-\int_{0}^{\infty}\int_{\Omega} N\partial_t\varphi-\int_{\Omega}\varphi(0,\cdot)N_0-d_1\int_{0}^{\infty}\int_{\Omega}N\Delta_x\varphi=\int_{0}^{\infty}\int_{\Omega}f(N,P)\varphi\\[0.1cm]
%&-\int_{0}^{\infty}\int_{\Omega} P\partial_t\psi-\int_{\Omega}\psi(0,\cdot)P_0-\int_{0}^{\infty}\int_{\Omega}\mu(N)P\Delta_x\psi=\int_{0}^{\infty}\int_{\Omega}g(N,P)\psi\\[0.2cm]
%&\forall\varphi,\psi \in \mathcal{C}^2_C([0,+\infty)\times\bar{\Omega}), \nabla_x \varphi \cdot n_{|_{\partial \Omega}}=0,\; \nabla_x \psi \cdot n_{|_{\partial \Omega}}=0.
%\end{align*}

Finally, $N$ satisfies the following extra regularity estimate: $N$ lies in $W^{1,2+\delta}([0,T]; L^{2+\delta} (\Omega)) \cap L^{2+\delta} ([0,T]; W^{2,2+\delta} (\Omega)) $  for all $T>0$ and some $\delta>0$.
% (and $P$ lies in $L^{2+\delta}([0,T]\times\Omega])$).
\medskip

Moreover, if $d=1$ or $d=2$ (remember that $d$ is the dimension of the domain) and if the initial data $N_{in}$, $P_{in}$ 
belong to $C^{0,\alpha}(\overline{\Omega})$ fore some $\alpha>0$, then all very-weak 
solutions of  (\ref{sfcrossS1}) -- (\ref{newint}) satisfy
$$N,\;P,\;\nabla_x N,\nabla_x P \in \mathcal{C}^{0,\alpha}([0,T]\times \bar{\Omega})$$
(for some $\alpha>0$), and
$$\partial_t N, \; \partial_t P,\; \partial_{x_ix_j}N,\; \partial_{x_ix_j}P \in L^p([0,T]\times \Omega) \quad \forall\; 1\leq p <\infty.$$
In other words, they are strong solutions.
\par 
Under the same assumptions on $d$, any couple of very-weak solutions $(N_1, P_1)$, $(N_2, P_2)$ with corresponding initial data $(N_{in, 1}, P_{in, 1})$,
$(N_{in, 2}, P_{in, 2})$ lying in $C^{0,\alpha}(\overline{\Omega})$ for some $\alpha>0$ satisfy the stability estimate
$$ ||N_1 - N_2||_{L^2([0,T] \times \Omega)} +
||P_1 - P_2||_{L^2([0,T] \times \Omega)} \le C_T\, (||N_{in,1} - N_{in,2}||_{L^2(\Omega)} + ||P_{in,1} - P_{in,2}||_{L^2(\Omega)} ) , $$
for some $C_T>0$ depending on $T$ (and on the data of the problem). This 
estimate ensures the uniqueness and stability of such very-weak solutions.
\par 
Finally, still under the same assumptions on $d$, and supposing that
$N_{in}$, $P_{in}$ 
belong to $C^{0,\alpha}(\overline{\Omega})$ fore some $\alpha>0$, and
$\inf {\hbox{ess }} N_{in}(x) >0$, the sequences $p_h^\var$ and  $p_s^\var$ converge a.e. towards $p_h$ and $p_s$.
% of  (\ref{sfcrossS1}) -- (\ref{newinit})
%and system \eqref{sfcrossS1}-\eqref{sfcrossS2}, together with homogeneous Neumann boundary conditions and initial data {\tobechecked $N_{in},\;P_{in} \in \mathcal{C}^{1,\alpha}(\bar{\Omega})$ %compatible with h.N.b.c.
%}, has a stable and unique solution for any $T>0$.
\end{thm}

%{\tobechecked
%\begin{prop}\label{thmsimplebis}
%Under the same assumptions as in Thm. \ref{thmsimple}, and under the extra assumptions that the dimension is $d=1$ or $d=2$, and that 

% Moreover  $p_h, p_s$ lie in $ L^{1} ([0,T]; W^{1,1} (\Omega))$ for all $T>0$.
%\end{prop}
%}
We next state the corresponding theorem in the case when $\xi>0$: 
\begin{thm}\label{thmadvanced}
Let $\Omega$ be a smooth domain of $\R^d$ (for some dimension $d \in \N - \{0\}$),  $d_1,\;d_2,\;d_3>0$ be diffusion rates, $r_0, \alpha, \xi, \tilde{\gamma}, \Gamma, \mu, \xi>0$ and $\eta\geq 0$ be parameters, and $N_{in}:=N_{in}(x) \ge 0$, $p_{h, in} := p_{h, in}(x) \ge 0$, $p_{s, in} := p_{s, in}(x) \ge 0$ be nonnegative initial data respectively in $L^{\infty}(\Omega)$, $L^{2+\delta}(\Omega)$, and $L^{2+\delta}(\Omega)$
for some $\delta>0$. We assume moreover that $\inf {\hbox{ess }} N_{in}(x) >0$.

Then for each $\var>0$, there exists a unique global classical (for $t>0$) solution ($N^\var$, $p_h^\var$, $p_s^\var$) of system \eqref{eqN1bis} -- \eqref{NBC}  (with the initial data defined above).

Moreover, when $\var \to 0$, one can extract from $N^\var$ a subsequence which is bounded in $L^{\infty}([0,T] \times\Omega)$ for all $T>0$ and converges a.e. towards a function $N\ge 0$ lying in  $L^{\infty}([0,T] \times\Omega)$, and from $p_s^\var$ (resp. $p_h^\var$) a subsequence which converges (strongly) in $L^{2+\delta}([0,T] \times\Omega)$ towards a function $p_s\ge 0$ (resp. $p_h\ge 0$) lying in  $L^{2+\delta}([0,T] \times\Omega)$ for all $T>0$ and some $\delta>0$.
% and from $p_h^\var$ a subsequence which converges (strongly) towards a function $p_h\ge 0$ lying in $L^{2+\delta}([0,T] \times\Omega)$ in $L^{2+\delta}([0,T] \times\Omega)$ for all $T>0$ and some $\delta>0$.

Moreover, $N$, $p_s$ and $p_h$ are very-weak solutions of the reaction-cross diffusion system
\begin{align}
\pa_t& N - d_1 \Delta_x N= r_0 \left(1 - \eta N\right) N - \frac{\alpha N p_s}{1 + \xi p_s}, \label{eqF1bis}\\
\pa_t& (p_s + p_h)  -  \Delta_x (d_2 p_s + d_3 p_h) =  \Gamma p_h - \mu (p_h + p_s), \label{eqF2bis}\\[0.2cm]
&\hspace{2.5cm}\frac{ \alpha N p_s}{1 + \xi p_s} = \tilde{\gamma} p_h ,\label{eqF3bis} 
\end{align}
with Neumann boundary conditions (\ref{NBC2})
% $\hat{n}(x)\cdot \nabla_x N=0,\quad \hat{n}(x)\cdot \nabla_x p_s=0, \quad \hat{n}(x)\cdot \nabla_x p_h=0,$
and with initial data (\ref{init}) in the following sense:
% $N(0,x) = N_{in}(x)$ and $(p_s + p_h)(0,x) = p_{s,in}(x)+p_{h,in}(x)$. 
identity (\ref{eqF3bis}) holds a.e., and for all
 $\varphi,\psi \in \mathcal{C}^2_c([0,+\infty)\times\bar{\Omega})$
 such that $ \nabla_x \varphi \cdot n_{|_{\partial \Omega}}=0,\; \nabla_x \psi \cdot n_{|_{\partial \Omega}}=0$,  
%\begin{multline*}
\begin{equation}\label{vw1bis}
-\int_{0}^{\infty}\int_{\Omega} N\partial_t\varphi-\int_{\Omega}\varphi(0,\cdot)N_{in}-d_1\int_{0}^{\infty}\int_{\Omega}N\Delta_x\varphi%\\
=\int_{0}^{\infty}\int_{\Omega} \bigg(r_0 (1-\eta N) N - \frac{\alpha N p_s}{1 +\xi\,p_s} \bigg)\, \varphi,
\end{equation}
%\end{multline*}
%\begin{multline*}
\begin{equation}\label{vw2bis}
-\int_{0}^{\infty}\int_{\Omega} (p_s + p_h)\partial_t\psi
-\int_{\Omega}\psi(0,\cdot)(p_{s,in}+p_{h,in})%\\
-\int_{0}^{\infty}\int_{\Omega}(d_2 p_s + d_3 p_h)\Delta_x\psi
=\int_{0}^{\infty}\int_{\Omega}(\Gamma p_h - \mu (p_h + p_s))\psi .
 \end{equation}
%\end{multline*}
%$$\forall\varphi,\psi \in \mathcal{C}^2_C([0,+\infty)\times\bar{\Omega}), \nabla_x \varphi \cdot n_{|_{\partial \Omega}}=0,\; \nabla_x \psi \cdot n_{|_{\partial \Omega}}=0.$$

Note that the reaction-cross diffusion system (\ref{eqF1bis}) -- (\ref{eqF3bis}) can be rewritten in the simpler form (\ref{sysBdAPDElim}), (\ref{suppl}), cf. computations of Subsection \ref{limsysder}.

Finally, $N$, $p_h$, $p_s$ satisfy the following extra regularity estimate: $N$ lies in $W^{1,p}([0,T]; L^{p} (\Omega))$ and 
$L^p ([0,T]; W^{2,p} (\Omega)) $  for all $T>0$ and all $p\in[1,+\infty[$, $p_h$ lies in $ L^{2} ([0,T]; H^{1} (\Omega)) $, and
$p_s$ lies in $ L^{1} ([0,T]; W^{1,1} (\Omega)) $.
\end{thm}

%{\tobechecked In this case, results on stability and uniqueness of the solution of \eqref{eqF1bis}-\eqref{eqF3bis} seem to be very difficult to obtain, due to the form of the cross-diffusion term, which involves $N,\;p_s,\;p_h$.}

\subsection{Study of the Turing instability}\label{TIa intro}
In Section \ref{SectionTIa} of this paper, we study the Turing instability regions associated to systems \eqref{sysBdAPDElim} and \eqref{sysBdAPDElin}. In order to do so, we first perform an adimensionalization, which enables to keep only a small number of parameters in the equations. Then we make explicit the condition on the parameters which leads to the existence of an homogeneous coexistence equilibrium for \eqref{sysBdAPDElim} and \eqref{sysBdAPDElin}. We also perform a linear stability analysis of this equilibrium (when it exists) at the level of ODEs (that is, w.r.t. homogeneous perturbations), and at the level of PDEs. Thus, we show that the Turing instability region (in terms of parameters) is nonempty, as expected, for both systems \eqref{sysBdAPDElim} and \eqref{sysBdAPDElin}. Finally, we compare the size of these regions. The main point is the fact that the Turing instability region associated to system \eqref{sysBdAPDElim} is always strictly included in the Turing instability region of system \eqref{sysBdAPDElin}. 

As a consequence, the use of reaction-diffusion systems for predator-prey interactions of Beddington-DeAngelis like in which standard diffusion is simply added to the reaction terms may lead to an overestimate of the possibility of appearance of patterns (at least in the case when the Beddington-DeAngelis functional response is a consequence of the interactions of searching and handling predators).

It is worth to mention that in some instances, the introduction of cross-diffusion terms instead of standard (linear) diffusion terms leads exactly to the opposite result, that is, the increase of the set of parameter values in which patterns develop, or even the appearance of patterns when none were observed with a linear diffusion \cite{Tulumello2014, iida}.
%%%%%%%%%%%%%%%%%%%%%%%%%%%%%%%%%%%%%%%%%%%%%%%%%%%%%%%%%%%%%%%%%%%%%%%%%%%%%%%%%%%%%%%%%%%%%%%%%%%%%%%%%%%%%%%%%%%%%%%%%%%%%%%%%%%%%%

\section{Rigorous results for the passages to the limit in microscopic models}\label{RigRes}

We systematically denote in this section by $C_T>0$ a constant which may depend on $T$ and on the parameters and initial data of the considered systems. \medskip

We start with the
\medskip

{\bf{Proof of Thm. \ref{thmsimple}}}\\
We consider system \eqref{eqN1bis} -- \eqref{eqN3bis} with $\xi=0$.
For a given $\var>0$, we  first recall that the existence of global in time solutions (for which $N^\var,\;p^\var_s,\;p^\var_h$ are nonnegative) for this system is classical (cf. \cite{D_milan} for example).\medskip 

Then, we observe that the r.h.s. of eq. \eqref{eqN1bis} is bounded above by $r_0N^\var$. Consequently, for each $T>0$, there exists $C_T>0$ such that 
 \begin{equation}\label{maxp}
  \sup_{\var>0} || N^\var ||_{L^{\infty}([0,T]\times \Omega)} \le C_T,
  \end{equation}
thanks to standard properties of the heat equation.
% Note that here, $C_T$ can in fact be bounded by $||N_{in}||_{L^\infty}e^{r_0T}$.
 As a consequence, there exists $N \in  L^{\infty}([0,T]\times \Omega)$ and a subsequence 
 (still denoted by $N^\var$) such that $N^\var \rightharpoonup N$ in $L^{\infty}([0,T]\times \Omega)$ weak $*$.\medskip
 
Adding the equations for $p_h^\var$ and $p_s^\var$, we end up with
\begin{equation}\label{eqsum1}
\pa_t P^\var  - \Delta_x (A^\var P^\var) =  \Gamma p_h^\var - \mu P^\var \le \Gamma P^\var,
\end{equation}
with 
\begin{equation}\label{plu}
 P^\var :=  p_h^\var  + p_s^\var, \qquad A^\var := \frac{ d_2 p_h^\var  + d_3 p_s^\var}{ p_h^\var  + p_s^\var}, 
 \end{equation}
so that thanks to a classical duality lemma (cf. \cite{DFPV} and the older reference \cite{PS}), 
$$ \sup_{\var>0} || P^\var ||_{L^{2}([0,T]\times \Omega)} \le C_T . $$
A refined version of the same lemma (cf. for example \cite{CDF} or \cite{BDF}) yields in fact the better estimate
\begin{equation}\label{deuxd} \sup_{\var>0} || P^\var ||_{L^{2+\delta}([0,T]\times \Omega)} \le C_T , 
 \end{equation}
for some $\delta>0$.

As a consequence, there exist $p_h, p_s \in  L^{2+\delta}([0,T]\times \Omega)$ and  subsequences (still denoted by $p_h^\var, p_s^\var$) such that $p_s^\var \rightharpoonup p_s$, $p_h^\var \rightharpoonup p_h$ in $L^{2+\delta}([0,T]\times \Omega)$ weak (for some $\delta>0$). \medskip 

Observing that $\pa_t N^\var - d_1 \Delta_x N^\var$ is bounded in  $L^{2+\delta}([0,T]\times \Omega)$ for all $T>0$ and some $\delta>0$,
we get thanks to the maximal regularity estimates for the heat kernel 
that
\begin{equation}\label{maxreg}
  \sup_{\var>0} || \pa_t N^\var ||_{L^{2+\delta}([0,T]\times \Omega)} \le C_T,  \qquad  \sup_{\var>0} || \pa_{x_i x_j} N^\var ||_{L^{2+\delta}([0,T]\times \Omega)} \le C_T. 
 \end{equation} 

We deduce from this estimate that the sequence $N^{\var}$ is strongly compact in $L^2([0,T]\times \Omega)$, so that (up to an extra extraction) $N^\var$ converges a.e. towards $N$. We also see that $N$ lies in $W^{1,2+\delta}([0,T]; L^{2+\delta} (\Omega)) \cap L^{2+\delta} ([0,T]; W^{2,2+\delta} (\Omega))$.\medskip

Using the bound (\ref{deuxd}), we end up with the convergence $N^\var p_s^\var \rightharpoonup N p_s$ in $L^{2+\delta}([0,T]\times \Omega)$ weak (for all $T>0$ and some $\delta>0$).\medskip

Passing to the limit in the equations \eqref{eqN1bis} and \eqref{eqsum1} in the sense of distributions, we end up with the equations \eqref{eqF1} and \eqref{eqF2}. More precisely, passing to the limit in the very-weak formulation  of equations \eqref{eqN1bis} and \eqref{eqsum1}, we get the very-weak formulation (\ref{vw1}), (\ref{vw2}) of the above system. %together with the homogeneous Neumann boundary condition on $N$ and $P$, and the initial data for $N$ and $P$ (that is $P(0,x)=p_{in}^s(x)+p_{in}^h(x)$).\medskip 

Observing that
$$ p_h^\var - N^\var p_s^\var = \var (\pa_t p_s^\var - d_2 \Delta_x p_s^\var) - \var (p_h^\var - \mu p_s^\var), $$
and passing to the limit in the sense of distributions in this statement, we get identity \eqref{eqF3}, which concludes the first part of the proof of Thm. \ref{thmsimple}.
\medskip

Now, we want to prove uniqueness and stability of very-weak solutions. In the rest of the proof of Thm. \ref{thmsimple}, we assume that $d=1$ or $d=2$. 
%\emph{Smoothness}\\
Consider now the equivalent form \eqref{sfcrossS1}-\eqref{sfcrossS2}
of system \eqref{eqF1}-\eqref{eqF3},
% together with homogeneous Neumann boundary conditions and initial data 
written  as
%in the following strong formulation
\begin{align}
\partial_t N - d_1 \Delta_x N = \phi(N) -\psi(N) P,&  \label{SFu1}\\[0.1cm]
\partial_t P - \Delta_x \left( \mu(N) P\right) = \nu(N) P,& \label{SFu2}
\end{align}
where $\phi,\;\psi, \; \mu, \nu: \mathbb{R}_+\to \mathbb{R}$ are  
defined by
$$\phi(N) := r_0(1-\eta N)N,\quad \psi(N) :=\dfrac{\tilde{\gamma}\alpha N}{\alpha N+\tilde{\gamma}},$$
$$\mu(N) :=\dfrac{d_2\tilde{\gamma}+d_3\alpha N}{\alpha N +\tilde{\gamma}},\quad \nu(N) :=\dfrac{\Gamma \alpha N}{\alpha N + \tilde{\gamma}}-\mu . $$
Then $\phi,\psi, \mu, \nu \in \mathcal{C}^2(\mathbb{R}^+)$, and $\mu \geq \mu_0>0$, with $\mu_0 := \inf(d_2, d_3)$.

Thanks to estimates (\ref{maxp}), (\ref{deuxd}) and (\ref{maxreg}),
we know  that $N \in L^{\infty} ([0,T]\times\Omega)$ and that $\partial_t N, \;  \partial_{x_ix_j} N,\; P \in L^{2+\delta} ([0,T]\times\Omega)$ for some $\delta>0$. 
%always $[0,T]\times\Omega$, not uniform in $t$.\\

%\medskip
%\noindent
%It can be proven that $N,\;P$ have much more regularity than hypothesis. In particular, in dimension 1 or 2
%$$N,\;P,\;\nabla_x N,\nabla_x P \in \mathcal{C}^{0,\alpha}([0,T]\times \bar{\Omega})$$
%for some $\alpha>0$, and
%$$\partial_t N, \; \partial_t P,\; \partial_{x_ix_j}N,; \partial_{x_ix_j}P \in L^p([0,T]\times \Omega) \quad \forall\; 1\leq p <\infty.$$
%Indeed,
%Maximal regularity for the heat equation for $N$\\
%I will use that $||\phi||_{L^\infty}\leq c,\; ||\psi||_{\infty}\leq c$.\\
%$\phi(N)-\psi(N)P \in L^{2+0}_{t,x},$ then
%$$\partial_t N \in L^{2+0}_{t,x},\quad \partial_{x_ix_j}N \in L^{2+0}_{t,x} \; \forall i,j$$
By interpolation, $\partial_{x_i} N \in L^{4+\delta}([0,T]\times\Omega)$ for any $i = 1,..,d$ and some $\delta>0$.
\par 
Computing
\begin{equation} \label{calmu}
\Delta_x \mu(N)=\nabla_x\cdot (\mu'(N)\nabla_x N)
             = \mu''(N)\Delta_x N +\mu'(N) |\nabla_x(N)|^2,
\end{equation}
we see that since $\mu\in \mathcal{C}^2(\mathbb{R}^+)$,  $\Delta_x N, |\nabla_x N|^2 \in L^{2+\delta}$ (for some $\delta>0$), then $\Delta_x \mu(N)\in L^{2+\delta}$ (for some $\delta>0$).\\
%now we need to know something on $P$, but it doesn't satisfy the heat equation!

%\subsection{(2st step)}
\medskip
%$$\partial_t P-\mu(N)\Delta_x P-2\nabla_x\mu(N)\nabla P - P\Delta_x\mu(N)=\nu(N)P.$$
Consider now any real number $q>0$ and compute (expanding 
$\Delta_x (\mu(N)\,P)$ and then performing integrations by parts):
% write the estimation with $q$ (use the equation, integrate by parts with the Stoke's formula).
\begin{align*}
&\partial_t \int_{\Omega}{\dfrac{P^{1+q}}{1+q}}=\int_{\Omega}{P^q\partial_t P}\\
&=\int_{\Omega}{\mu(N)P^q\Delta_xP}+2\int_{\Omega}{P^q\nabla_x\mu(N)\nabla_xP}+\int_{\Omega}{P^{1+q}\Delta_x\mu(N)}+\int_{\Omega}{\nu(N)P^{1+q}}\\
&=-\int_{\Omega}{\nabla_x (\mu(N) P^q)\cdot \nabla_x P} + 
2\int_{\Omega}{\nabla_x\mu(N)\cdot \nabla_x\left( \dfrac{P^{1+q}}{1+q}\right)}+\int_{\Omega}{P^{1+q}(\Delta_x\mu(N)+\nu(N))}\\
&=-\int_{\Omega}{\mu(N) q P^{q-1}|\nabla_x P|^2} + \int_{\Omega}{\nabla_x\mu(N)\cdot \nabla_x\left( \dfrac{P^{1+q}}{1+q}\right)}+\int_{\Omega}{P^{1+q}(\Delta_x\mu(N)+\nu(N))}\\
&=-q\int_{\Omega}{\mu(N) | P^{\frac{q-1}{2}}\nabla_x P|^2} -\int_{\Omega}{\Delta_x\mu(N)\dfrac{P^{1+q}}{1+q}} +\int_{\Omega}{P^{1+q}(\Delta_x\mu(N)+\nu(N))}\\
&=-q\int_{\Omega}{\mu(N)\left|\dfrac{\nabla_{x}(P^{\frac{q+1}{2}})}{\frac{q+1}{2}}\right|^2}+\int_{\Omega}{P^{1+q}\left(\dfrac{q}{1+q}\Delta_x\mu(N)+\nu(N)\right)} .\\
\end{align*}
Integrating  between $0$ and $T$, we end up with
%\begin{multline*}
\begin{equation}\label{esticomp}
\int_{\Omega}\dfrac{P^{1+q}(T,\cdot)}{1+q}+\dfrac{4q}{\left(q+1\right)^2}\int_{0}^T\int_{\Omega}{\mu(N)|\nabla_x(P^\frac{1+q}{2})|^2}=%\\
\int_{0}^{T}\int_{\Omega}{P^{1+q} \bigg( \dfrac{q}{1+q}\, \Delta_x\mu(N)+\nu(N) \bigg)}+\int_{\Omega}{\dfrac{P^{1+q}(0,\cdot)}{1+q}} .
\end{equation}
%\end{multline*}
%We have to bound the first term of the r.h.s. in order to get a good estimate on $P$. Note that it can be $P_{in}\in L^{1+q}$, (later $L^{q-1}$, $q$ to be discussed because it is not chosen).

%I will do an introduction.\\
Suppose now that $P\in L^a([0,T] \times\Omega)$ for some $a>2$.
Then $\phi(N)-\psi(N)P\in L^a([0,T] \times\Omega)$ and thanks to the maximal regularity estimates for the heat equation, we see that
%$$\phi(N)-\psi(N)P\in L^a_{t,x}$$
%as previously done replacing $2$ with $a$ ($a \neq \infty$, thanks to maximal regularity), we have
\begin{equation}\label{usef}
\partial_t N \in L^{a}([0,T] \times\Omega),\quad \partial_{x_ix_j}N \in L^{a}([0,T] \times\Omega).
\end{equation}
By interpolation, we see that 
% between $L^2$ and $L^\infty$,
 $\partial_{x_i}N\in L^{2a}([0,T] \times\Omega)$ and finally 
$$\Delta_x\mu(N)=\mu''(N)\Delta_x(N)+\mu'(N)|\nabla_xN|^2 \in L^{a}([0,T] \times\Omega).$$
Using estimate (\ref{esticomp}) for $q:= a-2$, we end up with
%We take $a-1=q-2$.
%\begin{multline*}
$$
\int_{\Omega}\dfrac{P^{a-1}(T, \cdot)}{a-1}+\dfrac{4(a-2)}{(a-1)^2}\int_{0}^T\int_{\Omega}{\mu(N)|\nabla_x(P^\frac{a-1}{2})|^2}=%\\
\int_{0}^{T}\int_{\Omega}{P^{a-1}\left(\frac{a-2}{a-1}\Delta_x\mu(N)+\nu(N)\right)}+\int_{\Omega}{\dfrac{P(0,\cdot)^{a-1}}{a-1}} .
$$
%\end{multline*}
Remembering that
$$P^{a-1}\in L^\frac{a}{a-1}([0,T] \times\Omega),\quad \Delta_x\mu(N) \in L^a([0,T] \times\Omega),\quad \nu(N)\in L^\infty([0,T] \times\Omega) ,$$%\quad P_{in}^{a-1}\in L^{a-1},
%then
%$$\Delta_x\mu(N)+\nu(N)\in L^a_{t,x}.$$
and observing that $L^a$ and $L^\frac{a}{a-1}$ are in duality, we see that
$$P^{a-1}\left(\Delta_x\mu(N)+\nu(N)\right)\in L^1([0,T] \times\Omega),$$
so that
%$$\int_{\Omega} \frac{P^{a-1}}{a-1}(T)\leq C_T$$
\begin{equation}
\int_{\Omega} \frac{P^{a-1}(T,\cdot)}{a-1} + \mu_0\int_{0}^T\int_{\Omega}{|\nabla_x(P^\frac{a-1}{2})|^2}
%\leq\int_{0}^T\int_{\Omega}{\mu(N)|\nabla_x(P^\frac{a-1}{2})|^2}
\leq C_T .
\label{eq:star1}
\end{equation}
As a consequence $P\in L^\infty([0,T]; L^{a-1}(\Omega))$ and $P^\frac{a-1}{2}\in L^2([0,T]; H^1(\Omega))$. 
%(it has a gradient in $L^2$). Now we 
Using Sobolev estimates, we see that
\smallskip

 if $d=1$,  $P^\frac{a-1}{2}\in L^2([0,T]; L^{\infty}(\Omega))\quad \Rightarrow \quad P\in L^{a-1}([0,T]; L^\infty(\Omega))$;
 \smallskip 
 
 if $d=2$,  $P^\frac{a-1}{2}\in L^2([0,T]; L^{q}(\Omega))$ for all
 $q \in [1, \infty[ \quad \Rightarrow \quad  P\in L^{a-1}([0,T]; L^q(\Omega))$  for all
 $q \in [1, \infty[$.
 %   P\in L^{a-1}_t(L^{\infty-0}_x)$.
%\item[d=3] $P^\frac{a-1}{2}\in L^2_t(L^{6}_x)\quad \rightarrow \quad P\in L^{a-1}_t(L^{3(a-1)}_x)$.
%\end{description}
\smallskip

\noindent
By interpolation with $L^\infty([0,T]; L^{a-1}(\Omega))$,
\smallskip 

%\begin{description}
 if $d=1$,   $P\in L^{2(a-1)}([0,T] \times \Omega)$;
 \smallskip 
 
 if $d=2$,  $P\in L^{2(a-1)-\delta} ([0,T] \times \Omega) $ for any 
 $\delta>0$.
%\item[d=3] In this case it is not enough.%$P \in$ something?.
%\end{description}
%{\tobechecked For now, we keep $d=1,2$.}
\smallskip

\noindent
We define the sequence $(a_n)_{n\in\N}$ by 
%$$\begin{cases}
$a_0>2$; $a_{n+1}=2(a_n-1)$, and observe that 
%\end {cases}$$
 $a_n\to\infty$.
 % ({\tobechecked $a_n=2^n(\dots)+2$}). 
\par
Starting from estimate (\ref{deuxd}) and proceeding 
by induction, we see that $P\in L^{2(a_n-1)-\delta} ([0,T] \times \Omega) $ for all $n\in\N$,  and also that
\begin{equation}
\partial_t{N},\; \partial_{x_ix_j}N\in L^{q}([0,T] \times \Omega)
\label{eq:starsmooth}
\end{equation}
for all $q \in [1, \infty[$ thanks to (\ref{usef}).
%Starting from $P\in L^2$ we have $P\in L^{???}$: in fact $P\in L^\infty_t(L^{\infty-0})_x,\; \nabla_x P \in L^2_{t,x}$, taking $a=3$ in \eqref{eq:star1}. 
%Questions: is it possible to have better results for $P$ and $N$? For the uniqueness, we need $\Delta_x \mu(N) \in L^\infty,\; \nabla_x P\in L^\infty$. What we have for now is not enough!

%\subsection{(3rd step)}
%We have that 
%$$\nabla_x(\phi(N)-\phi(N)P)=\phi'(N)\nabla_xN-\psi'(N)\nabla_x N P-\psi(N)\nabla_x P.$$
%Since
%$$\phi'(N)\in L^\infty,\; \nabla_x N \in L^{4+0},\; \psi'(N) \in L^\infty,\; P\in L^{\infty-0},\; \psi(N)\nabla_xP\in L^2_{t,x},$$
%we have that
%$$\phi'(N)\nabla_xN \in L^{4+0}_{t,x},\; \psi'(N)\nabla_x N P\in L^{4+0}_{t,x},$$
%then the r.h.s. belongs to $L^2_{t,x}$. As a consequence, the r.h.s. of the first equation has gradient in $L^2$. Thanks to the maximal regularity (to be checked) $\partial_t N-d_1\Delta_x N$ has gradient in $L^2$ and $\Delta_x \mu(N)\in L^\infty_{t,x}$. 
%Thanks to a result in \cite{ladyzhenskaia1988} (Chapter IV, theorem 9.1 and the corollary).
%$$\partial_t P-\mu(N)\Delta_x P-2\nabla_x\mu(N)\nabla_xP-(\Delta_x\mu(N)-\nu(N))P=0.$$
Thanks to identity (\ref{calmu}) and the properties of $\mu$, 
we know therefore that $\nabla_x\mu(N)\in L^{\infty}([0,T] \times \Omega)$, and $\Delta_x\mu(N)+\nu(N) \in L^{q}([0,T] \times \Omega)$ for 
all $q \in [1, \infty[$.
% thanks to \eqref{eq:starsmooth}. If $\mu(N)\geq c_0>0$ and continuous in $t,\;x$ %(these are hypothesis) 
%then $\partial_t P$ and $\partial_{x_ix_j}P \in \L_{t,x}^{\infty-0}$. We need maximal regularity for the heat equation
%\begin{equation}
%\phi(N)-\psi(N)P \in L^{\infty-0}_{t,x}.
%\label{eq:0}
%\end{equation}
%Then both $\partial_t N$ and $\partial_{x_ix_j}N$ lie in $L^{\infty-0}_{t,x}$. By Sobolev embedding, $N$ is continuous in $t,\;x$ and $\nabla_x N\in L^\infty_{t,x}$. 
%More in detail, we have that $\phi(N),\; \psi(N) \in L^{\infty}_{t,x}, \; P \in L^{\infty-0}_{t,x}$. Then $\phi(N)-\psi(N)P \in L^{\infty-0}_{t,x}$ (this is \eqref{eq:0}). For maximal regularity, $\partial_t N,\;\partial_{x_ix_j}\in L^{\infty-0}_{t,x} \; \forall i,j$. Thanks to properties of the heat equation ({\tobechecked add some references to Laurent's papers}), we have that $\partial_{x_i}\in L^{\infty}_{t,x}\; \forall i$. Thus, 
%Then
%$$\nabla_x\mu(N)=\mu'(N)\nabla_x(N)\nabla_x N \in L^{\infty-0}_{t,x}$$
%and
%$$\Delta_x\mu(N)=\mu''(N)\Delta_x N+\mu'(N)|\nabla_xN|^2 \quad \in L^{\infty-0}_{t,x},$$
%because $\mu''(N),\;\mu'(N),\;|\nabla_xN|^2\; \in L^{\infty}_{t,x}$ and $\Delta_x N\in L^{\infty-0}_{t,x}$. 
%Then, $\partial_t N \in L^{\infty-0}_{t,x}$ and $\partial_{x_i}\in L^{\infty}_{t,x}$ imply that $N$ is continuous by Sobolev embedding.
Expanding eq. (\ref{SFu2}) as 
$$ \pa_t P - \mu(N)\, \Delta_x P - 2\, \nabla_x \mu(N) \cdot \nabla_x P
= (\Delta_x \mu(N) + \nu(N))\, P, $$
and using
 Theorem 9.1 and its corollary in \cite{ladyzhenskaia1988} (see also \cite{DeTr}), we get for $P$ the estimates 
$$\partial_t P \in L^{q}([0,T] \times \Omega), \qquad \partial_{x_ix_j}P \in L^{q}([0,T] \times \Omega) ,$$
for 
all $q \in [1, \infty[$, $i,j \in \{1,..,d\}$.
\par 
Thanks once again to Sobolev embeddings, we also see that
%then for properties of the heat equation \cite{CDF}
$$P,\; \nabla_x P \in L^{\infty}([0,T] \times \Omega) \qquad  \forall i,j \in \{1,..,d\} . $$
% \; \forall i.$$
%{\tobechecked Is it enough now to have uniqueness? We need $\Delta_xN \in L^{\infty}_{t,x}$ in order to have $\Delta_x\mu(N) \in L^{\infty}_{t,x}$.} 
\medskip 

\noindent
%\textit{Stability and uniqueness}\\
%How to use these results to get stability (close to initial data forever) and uniqueness?
We now prove the statement about stability, still assuming that the dimension
is $d=1$ or $d=2$. 
Let $N_i,\; P_i,\; i=1,2$ be two different solutions of the same problem \eqref{SFu1} -- \eqref{SFu2}, both with homogeneous Neumann boundary condition, but with different initial data: %(this gives the stability)
\begin{align}
\partial_t N_i - d_1 \Delta_x N_i = \phi(N_i) -\psi(N_i) P_i,\\[0.1cm]
\partial_t P_i - \Delta_x \left( \mu(N_i) P_i\right) = \nu(N_i) P_i,\\[0.2cm]
\nabla_x N_i \cdot n_{|_{\partial \Omega}}=0,\; \nabla_x P_i \cdot n_{|_{\partial \Omega}}=0, \\[0.1cm]
N_i(0,\cdot)=N^i_{in},\; P_i(0,\cdot)=P^i_{in}.
\end{align}
We can write a first estimate for $N_1-N_2$. We compute
$$\partial_t (N_1-N_2) - d_1 \Delta_x (N_1-N_2) = \phi(N_1) -\phi(N_2)-\psi(N_1) P_1+\psi(N_2) P_2,$$
multiply by $(N_1-N_2)$ this formula,  and then integrate w.r.t. space (plus an integration by parts and the use of the Young's inequality); we get
\begin{multline} \label{trq}
\partial_t \int_{\Omega}{\dfrac{(N_1-N_2)^2}{2}}+ d_1\int_{\Omega}{|\nabla_x(N_1-N_2)|^2}=
\int_{\Omega}{(\phi(N_1) -\phi(N_2))(N_1-N_2)}\\
-\int_{\Omega}{P_1(\psi(N_1)-\psi(N_2))(N_1-N_2)}
+\int_{\Omega}{\psi(N_2)(P_2-P_1)(N_1-N_2)}\\
\leq 
||\phi'||_{L^{\infty}([0, \sup(||N_1||_{L^{\infty}}, 
||N_2||_{L^{\infty}} )] )}\, \int_{\Omega}{(N_1-N_2)^2}
+||\psi'||_\infty ||P_1||_\infty \int_{\Omega}{(N_1-N_2)^2}
+\dfrac{||\psi||_\infty}{2}\int_{\Omega}{(P_1-P_2)^2}\\
+\dfrac{||\psi||_\infty}{2}\int_{\Omega}{(N_1-N_2)^2}
\leq 
K_2 \left( \int_{\Omega}{(N_1-N_2)^2} +\int_{\Omega}{(P_1-P_2)^2}\right),
\end{multline}
where $K_2:=K_2(||\phi'||_{L^{\infty}([0, \sup(||N_1||_{L^{\infty}}, 
||N_2||_{L^{\infty}} )] )},||\psi'||_\infty||P_1||_{L^{\infty}},||\psi||_\infty)$. 

We also compute
$$\partial_t(P_1-P_2)-\Delta_x(\mu(N_1)P_1-\mu(N_2)P_2)=\nu(N_1)P_1-\nu(N_2)P_2,$$
which can be multiplied by $(P_1-P_2)$ and then integrated in space. We get
\begin{multline*}
\partial_t \int_{\Omega}{\dfrac{(P_1-P_2)^2}{2}}+\int_{\Omega}{\nabla_x(P_1-P_2)\cdot\nabla_x [ \mu(N_1)(P_1-P_2)+P_2(\mu(N_1)-\mu(N_2)) ] }\\=\int_{\Omega}{(\nu(N_1)P_1-\nu(N_2)P_2)(P_1-P_2)}.
\end{multline*}
It can be rewritten as 
\begin{multline*}
\partial_t \int_{\Omega}{\dfrac{(P_1-P_2)^2}{2}}
+\int_{\Omega}{\mu(N_1)|\nabla_x(P_1-P_2)|^2}=\\
-\int_{\Omega}{(P_1-P_2)\nabla_x(P_1-P_2)\cdot\nabla_x\mu(N_1)}
-\int_{\Omega}{P_2\nabla_x(P_1-P_2)\cdot\nabla_x(\mu(N_1)-\mu(N_2))}\\
-\int_{\Omega}{(\mu(N_1)-\mu(N_2))\nabla_x(P_1-P_2)\cdot\nabla_xP_2}
+\int_{\Omega}{\nu(N_1)(P_1-P_2)^2}
+\int_{\Omega}{P_2(\nu(N_1)-\nu(N_2))(P_1-P_2)}.
\end{multline*}
We observe that (for any $\var>0$)
\begin{multline*}
-\int_{\Omega}{(P_1-P_2)\nabla_x(P_1-P_2)\cdot\nabla_x [\mu(N_1)] }\\
\leq
||\nabla_x [\mu(N_1)] ||_{L^\infty}\dfrac{\varepsilon}{2}\int_{\Omega}{|\nabla_x(P_1-P_2)|^2}
+||\nabla_x [\mu(N_1)] ||_{L^\infty}\dfrac{1}{2\varepsilon}\int_{\Omega}{(P_1-P_2)^2},
\end{multline*}

\begin{multline*}
 - \int_{\Omega}{P_2\nabla_x(P_1-P_2)\cdot\nabla_x(\mu(N_1)-\mu(N_2))} \\
\leq
||P_2||_{L^\infty}\varepsilon\int_{\Omega}{|\nabla_x(P_1-P_2)|^2}
+||P_2||_{L^\infty}\dfrac{||\mu'||^2_{\infty}}{2\varepsilon}\int_{\Omega}{|\nabla_x(N_1-N_2)|^2}\\
+||P_2||_{L^\infty}\dfrac{||\mu''||^2_{\infty}}{2\varepsilon}||\nabla_xN_2||^2_{L^\infty}\int_{\Omega}{(N_1-N_2)^2},
\end{multline*}

\begin{multline*}
- \int_{\Omega}{(\mu(N_1)-\mu(N_2))\nabla_x(P_1-P_2)\cdot\nabla_xP_2}
\leq
||\nabla_xP_2||_{L^\infty}\dfrac{\varepsilon}{2}\int_{\Omega}{|\nabla_x(P_1-P_2)|^2}\\
+||\nabla_xP_2||_{L^\infty}\dfrac{||\mu'||^2_{L^\infty}}{2\varepsilon}\int_{\Omega}{(N_1-N_2)^2} ,
\end{multline*}

%\begin{multline*}
$\displaystyle{
\int_{\Omega}{\nu(N_1)(P_1-P_2)^2}
\leq
 ||\nu||_{\infty} \int_{\Omega}{(P_1-P_2)^2}
\dfrac{||P_2||_{L^\infty}}{2}\int_{\Omega}{(P_1-P_2)^2}%\\
+\dfrac{||P_2||_{L^\infty}}{2}||\nu'||^2_{L^\infty}\int_{\Omega}{(N_1-N_2)^2},
}
$\\[0.2cm]
%\end{multline*}

$\displaystyle{
\int_{\Omega}{P_2(\nu(N_1)-\nu(N_2))(P_1-P_2)}
\leq \dfrac{||P_2||_{L^\infty}}{2}\int_{\Omega}{(P_1-P_2)^2}%\\
+\dfrac{||P_2||_{L^\infty}}{2}||\nu'||^2_{L^\infty}\int_{\Omega}{(N_1-N_2)^2} .
}
$\\[0.2cm]

%and then
%\begin{multline*}
%\partial_t \int_{\Omega}{\dfrac{(P_1-P_2)^2}{2}}
%+\int_{\Omega}{\mu(N_1)|\nabla_x(P_1-P_2)|^2}=
%-\int_{\Omega}{(P_1-P_2)\nabla_x(P_1-P_2)\cdot\nabla_x\mu(N_1)}
%-\int_{\Omega}{P_2(P_1-P_2)\nabla_x(P_1-P_2)\cdot\nabla_x(\mu(N_1)-\mu(N_2))}\\
%-\int_{\Omega}{(\mu(N_1)-\mu(N_2))\nabla_x(P_1-P_2)\cdot\nabla_xP_2}+\int_{\Omega}{\nu(N_1)(P_1-P_2)^2}+\int_{\Omega}{P_2(\nu(N_1)-\nu(N_2))(P_1-P_2)}\\
%\leq 
%-\int_{\Omega}{(P_1-P_2)\nabla_x(P_1-P_2)\cdot\nabla_x\mu(N_1)}+||P_2||_{L^\infty}\epsilon\int_{\Omega}{|\nabla_x(P_1-P_2)|^2}+
%||P_2||_{L^\infty}\dfrac{||\mu'||_{\infty}}{2\epsilon}\int_{\Omega}{|\nabla_x N_1-\nabla_x N_2|^2}\\
%+||P_2||_{L^\infty}\dfrac{||\nabla_x N_2||_{\infty}}{2\epsilon}||\mu''||^2_{\infty}\int_{\Omega}{(N_1- N_2)^2}\\
%+||\nu||_{\infty}\int_{\Omega}{(P_1- P_2)^2}+\dfrac{||P_2||_{L^\infty}}{2}\int_{\Omega}{(P_1- P_2)^2}+\dfrac{||P_2||_{L^\infty}}{2}||\nu'||^2_{\infty}\int_{\Omega}{(N_1- N_2)^2}.
%\end{multline*}
%%When $\epsilon$ is small enough I can eliminate two terms. I have to bound the first term of the r.h.s., 
%\begin{multline*}
%-\int_{\Omega}{(P_1-P_2)\nabla_x(P_1-P_2)\cdot\nabla_x\mu(N_1)}\\
%\leq
%||\nabla_x\mu(N_1)||_{L^\infty}\dfrac{\epsilon}{2}\int_{\Omega}{|\nabla_x(P_1-P_2)|^2}
%+||\nabla_x\mu(N_1)||_{L^\infty}\dfrac{1}{2\epsilon}\int_{\Omega}{(P_1-P_2)^2}.
%\end{multline*}
%Then we have that
Using these inequalities, we get
\begin{multline}
\partial_t \int_{\Omega}{\dfrac{(P_1-P_2)^2}{2}}
+\int_{\Omega}{\mu(N_1)|\nabla_x(P_1-P_2)|^2}\\
\leq 
K_1\left(\dfrac 3 2 + \dfrac 1 {2\varepsilon}\right) \int_{\Omega}(P_1-P_2)^2 
+K_1\left(\dfrac 1 2 + \dfrac 1 {\varepsilon}\right) \int_{\Omega}(N_1-N_2)^2 \\
+K_1 2\varepsilon \int_{\Omega}{|\nabla_x (P_1-P_2)|^2}
+K_1 \dfrac{1}{2\varepsilon}\int_{\Omega}{|\nabla_x (N_1-N_2)|^2},
\end{multline}
where 
\begin{multline*}
K_1:=\max(||\nu||_{\infty}, ||\mu||_{\infty}, ||P_2||_{L^\infty}, ||\nabla_x\mu(N_1)||_{L^\infty}, ||P_2||_{L^\infty}||\nu'||^2_{\infty},
 ||\nabla_x P_2||_{L^\infty}||\mu'||^2_{\infty},\\
 ||P_2||_{L^\infty}||\mu''||^2_{\infty}||\nabla_x N_2||_{L^\infty}^2,
||P_2||_{L^\infty}||\mu'||^2_{\infty},||\nabla_xP_2||_{L^\infty}).
\end{multline*}
\medskip

Remembering (\ref{trq}), we end up with the system of differential inequalities:
%\begin{multline*}
$$
\partial_t \int_{\Omega}{\dfrac{(N_1-N_2)^2}{2}}+ d_1\int_{\Omega}{|\nabla_x(N_1-N_2)|^2}%\\
\leq 
K_2 \left( \int_{\Omega}{(N_1-N_2)^2} +\int_{\Omega}{(P_1-P_2)^2}\right),
$$
%\end{multline*}
and
\begin{multline*}
\partial_t \int_{\Omega}{\dfrac{(P_1-P_2)^2}{2}}+ \mu_0\int_{\Omega}{|\nabla_x(P_1-P_2)|^2}\\
\leq 
K_1\left(\left(\dfrac 3 2 + \dfrac 1 {2\varepsilon}\right) \int_{\Omega}(P_1-P_2)^2 
+\left(\dfrac 1 2 + \dfrac 1 {\varepsilon}\right) \int_{\Omega}(N_1-N_2)^2 
+2\varepsilon \int_{\Omega}{|\nabla_x (P_1-P_2)|^2}
+\dfrac{1}{2\varepsilon}\int_{\Omega}{|\nabla_x (N_1-N_2)|^2}\right),
\end{multline*}
which holds for all $\varepsilon>0$. Taking 
$\varepsilon := \dfrac{\mu_0}{2K_1}$,
 the second inequality becomes
\begin{multline*}
\partial_t \int_{\Omega}{\dfrac{(P_1-P_2)^2}{2}} + \mu_0\int_{\Omega}{|\nabla_x(P_1-P_2)|^2}
\leq 
\left(\dfrac {3K_1} 2 + \dfrac {K_1^2}{\mu_0}\right) \int_{\Omega}(P_1-P_2)^2 \\
+\left(\dfrac {K_1} 2 + \dfrac {2K_1^2}{\mu_0}\right) \int_{\Omega}(N_1-N_2)^2 
+ \mu_0 \int_{\Omega}{|\nabla_x (P_1-P_2)|^2}
+\dfrac{K_1^2}{\mu_0}\int_{\Omega}{|\nabla_x (N_1-N_2)|^2} .
\end{multline*}
We now consider a linear combination of the two inequalities: 
%(the first multiplied by $K_1^2/c_0$, the second by $d_1$), obtaining
\begin{multline*}
\partial_t \left( \dfrac{K_1^2}{2 \mu_0}\int_{\Omega}{(N_1-N_2)^2}+\frac{d_1}{2}\int_{\Omega}{(P_1-P_2)^2}\right)\\
\leq 
\left(\dfrac {K_1^2K_2} {\mu_0} + d_1K_1\left(\dfrac{1}{2}+\dfrac {2K_1} {\mu_0}\right)\right) \int_{\Omega}(N_1-N_2)^2 
+\left(\dfrac {K_1^2K_2} {\mu_0} + d_1K_1\left(\dfrac{3}{2}+\dfrac {K_1} {\mu_0}\right)\right)  \int_{\Omega}(P_1-P_2)^2.
\end{multline*}
Choosing 
$$C_1:=\min{\left(\dfrac{K_1^2}{2\mu_0},\frac{d_1}{2}\right)},\quad C_2 := \max{\left[\dfrac {K_1^2K_2} {\mu_0} + d_1K_1\left(\dfrac{1}{2}+\dfrac {2K_1} {\mu_0}\right),\dfrac {K_1^2K_2} {\mu_0} + d_1K_1\left(\dfrac{3}{2}+\dfrac {K_1} {\mu_0}\right)\right]},$$
we end up with
$$\partial_t \left( \int_{\Omega}{(N_1-N_2)^2}+\int_{\Omega}{(P_1-P_2)^2}\right)\\
\leq 
\dfrac{C_2}{C_1}\left( \int_{\Omega}{(N_1-N_2)^2}+\int_{\Omega}{(P_1-P_2)^2}\right).$$
%Integrating with respect to time, we have
%\begin{multline*}
%\int_{\Omega}{(N_1-N_2)^2}(T)+\int_{\Omega}{(P_1-P_2)^2}(T)
%\leq 
%\int_{\Omega}{(N_1-N_2)^2}(0)+\int_{\Omega}{(P_1-P_2)^2}(0)\\
%+\dfrac{C_2}{C_1}\left(\int_{0}^{T}\int_{\Omega}{(N_1-%N_2)^2}+\int_{0}^{T}\int_{\Omega}{(P_1-P_2)^2}\right),
%\end{multline*}
Using finally Gronwall's lemma, we get the statement of  stability and
therefore also of uniqueness in Thm. \ref{thmsimple}.
\medskip

We now prove the last statement of Thm. \ref{thmsimple}:
%{\bf{Proof of Prop. \ref{thmsimplebis}}}\\
We compute (for $q \in ]0,1[$) the derivative of the following nonnegative quantity:
\begin{multline}
\frac{d}{dt}\int_{\Omega} \bigg[\frac{(p_h^\var)^{q+1}}{q+1}+\bigg(\frac{\alpha}{\tilde{\gamma}}N^\var\bigg)^q\frac{(p_s^\var)^{q+1}}{q+1} \bigg]\\
= \int_{\Omega}  (p_h^\var)^{q} [ d_3\, \Delta_x p_h^\var 
- \frac1{\var} ( - \alpha\, N^\var  p_s^\var + \tilde{\gamma}\, p_h^\var )
- \mu\, p_h^\var ] \\
+  \int_{\Omega} \bigg( \frac{\alpha}{\tilde{\gamma}} N^\var \bigg)^q\,
(p_s^\var)^q\, [ d_2\, \Delta_x p_s^\var 
+  \frac1{\var} ( - \alpha\, N^\var  p_s^\var + \tilde{\gamma}\, p_h^\var  ) 
+ \Gamma\, p_h^\var - \mu\, p_s^\var ] \\
+ \bigg( \frac{\alpha}{\tilde{\gamma}} \bigg)^q
\,  \int_{\Omega} q\, (N^\var)^{q-1} \frac{(p_s^\var)^{q+1}}{q+1} \, \pa_t 
N^\var  \\
= - \underbrace{d_3 q \int_{\Omega} (p_h^\var)^{q-1}|\nabla_x p_h^\var|^2}_{\boxed{1}}
  - \underbrace{d_2 q \left(\dfrac{\alpha}{\tilde{\gamma}}\right)^q \int_{\Omega}  (N^\var)^q (p_s^\var)^{q-1} |\nabla_x p_s^\var|^2}_{\boxed{2}}  \\
  - \underbrace{\frac{\tilde{\gamma}}{\var} \int_{\Omega} \bigg(p_h^\var - \left(\dfrac{\alpha}{\tilde{\gamma}}\right) N^\var p_s^\var  \bigg)\, \left((p_h^\var)^q - \left(\dfrac{\alpha}{\tilde{\gamma}} N^\var  p_s^\var\right)^q \right)}_{\boxed{3}} \\
  - \underbrace{\mu \int_{\Omega} \bigg[ (p_h^\var)^{q+1} + \left(\dfrac{\alpha}{\tilde{\gamma}} N^\var\right)^q(p_s^\var)^{q+1} \bigg]}_{\boxed{4}} 
  + \underbrace{\Gamma \left(\dfrac{\alpha}{\tilde{\gamma}}\right)^q
   \int_{\Omega} (N^\var)^q (p_s^\var)^q p_h^\var}_{\boxed{5}} \\
  + \underbrace{\frac{q}{q+1} \left(\dfrac{\alpha}{\tilde{\gamma}}\right)^q \int_{\Omega} (p_s^\var)^{q+1} (N^\var)^{q-1} (\pa_t N^\var + d_2 \Delta_x N^\var)}_{\boxed{6}}\\
  - \underbrace{d_2 \frac{q(1-q)}{1+q} \left(\dfrac{\alpha}{\tilde{\gamma}}\right)^q \int_{\Omega} (p_s^\var)^{q+1} (N^\var)^{q-2} |\nabla_x N^\var|^2 }_{\boxed{7}}.\label{splitder}
\end{multline}
We observe that the terms \boxed{1}, \boxed{2}, \boxed{3}, \boxed{4}, \boxed{7} are all nonpositive. Remembering that $N^\var$ is bounded in $L^\infty$, and that $p_s^\var,\;p_h^\var$ are bounded in $L^{2+\delta}$ for some $\delta>0$ (see estimates (\ref{maxp}) and (\ref{deuxd})), we see that
$$\int_{0}^{T}{\int{(N^\var)^q(p_s^\var)^qp_h^\var}}\leq C_T \quad \textnormal{ (for all $q\in]0,1[$)}.$$
Remembering then that $\partial_t N^\var+d_2\Delta_x N^\var$ is bounded in $L^{2+\delta}$ for some $\delta>0$ (see estimate (\ref{maxreg})), we see that
$$\int_{0}^{T}{\int{(p_s^\var)^{q+1}\, (N^\var)^{q-1}\, |\partial_t N^\var+d_2\Delta_x N^\var|}} \leq C_T\quad \textnormal{ when } q\in]0,1[ \textnormal{ is small enough}.$$
As a consequence, integrating \eqref{splitder} on $[0,T]$, we see that (for $q\in]0,1[$ small enough) 
\begin{equation} 
\int_0^T \int_{\Omega} \bigg(p_h^\var - \frac{\alpha}{\tilde{\gamma}} N^\var p_s^\var  \bigg)\, \bigg((p_h^\var)^q - ( \frac{\alpha}{\tilde{\gamma}} N^\var
 p_s^\var)^q \bigg) \le C_T\var, \label{in1}
\end{equation} 
$$\int_0^T\int_{\Omega} (p_h^\var)^{q-1} |\nabla_x p_h^\var|^2 \le C_T, $$ 
$$ \int_0^T\int_{\Omega} (N^\var)^q  (p_s^\var)^{q-1} |\nabla_x p_s^\var|^2 \le C_T. $$
Observing that
\begin{equation} (\pa_t - d_1 \Delta_x) \ln N^\var=\dfrac{1}{N^\var}(\pa_t -d_1\Delta_x)N^\var+d_1\dfrac{|\nabla_xN^\var|^2}{(N^\var)^2} \ge r_0(1-\eta N^\var)- \alpha p_s^\var \geq (-\eta r_0 - \alpha p_s^\var) C_T, 
\end{equation}
we see that since $p_s^\var$ is bounded in $L^{2+\delta}([0,T]\times \Omega)$ for some $\delta>0$, in dimension $d=1$ or $d=2$, we obtain
that $N^\var$ is bounded below (by a strictly positive constant) on 
$[0,T] \times \Omega$ as soon as $\inf ess\, N_{in} >0$. Indeed, we recall that $(\pa_t -d_1\Delta_x)^{-1}$ acts as a convolution with a function lying in $L^{3-\var}$ (when $d=1$) or $L^{2-\var}$ (when $d=2$) for any $\var>0$ (cf. \cite{CDF}), so that thanks to the Young's inequality and the assumption that the initial datum $N^\var$ is essentially strictly positive, $\ln N^\var$ is bounded below (by a strictly positive constant). As a consequence, still for $q\in]0,1[$ small enough,
$$ \int_0^T\int_{\Omega} (p_s^\var)^{q-1}  |\nabla_x p_s^\var|^2 \le C_T . $$

Then, Cauchy-Schwarz inequality ensures that 
\begin{equation} 
\bigg( \int_0^T\int_{\Omega} |\nabla_x p_{s,h}^\var| \bigg)^2 \le  \bigg( \int_0^T\int_{\Omega}  (p_{s,h}^\var)^{q-1}  |\nabla_x p_{s,h}^\var|^2
\bigg)      \times  \, \bigg(\int_0^T\int_{\Omega} (p_{s,h}^\var)^{1-q}
\bigg) \le C_T. \label{inCS}
\end{equation}							

Using \eqref{eqsum1}, we see that $\pa_t P^\var$ is bounded in $L^{2}([0,T]; H^{-2}( \Omega))$, so that thanks to Aubin's lemma \cite{moussa}, $P^\var$ converges a.e. to $P$.
Then we use the elementary inequality (which holds for all $x,y\ge 0$ and some
constant $C>0$)
$$ (x-y)\, (x^q - y^q) \ge C\, (x^{\frac{1+q}2} - y^{\frac{1+q}2})^2 ,$$
and extract from (\ref{in1}) the estimate
\begin{equation} 
\int_0^T \int_{\Omega} \bigg( (p_h^\var)^{\frac{1+q}2}  -  \left(\frac{\alpha}{\tilde{\gamma}} N^\var p_s^\var\right)^{\frac{1+q}2}  \bigg)^2\,  \le C_T\var .
\end{equation} 
Using another elementary inequality (still holding for all $x,y\ge 0$, and 
$q>0$ small enough), namely
$$ |x-y| \le C\,  |x^{\frac{1+q}2} - y^{\frac{1+q}2}|
\times (x^{\frac{1-q}2} - y^{\frac{1-q}2}) ,$$
and Cauchy-Schwarz inequality, we see that
\begin{equation} 
\int_0^T \int_{\Omega} |p_h^\var  -  \left(\frac{\alpha}{\tilde{\gamma}}\right) N^\var p_s^\var|  \le C_T\,\sqrt{\var} \times
\bigg( \int_0^T \int_{\Omega} \bigg[ (p_h^\var)^{1-q}  +
  \left(\frac{\alpha}{\tilde{\gamma}} N^\var p_s^\var\right)^{1-q} \bigg]\,  \bigg)^{1/2}
\le   C_T\,\sqrt{\var}  .
\end{equation} 

 Then  $ p_h^\var - \frac{\alpha}{\tilde{\gamma}}\, N^\var p_s^\var \to 0$ a.e. (up to extraction of a subsequence). Remembering that $N^\var$ converges a.e. to $N$ and $P^\var$ converges a.e. to $P$, we see that $p_h^\var$ converges a.e. to $p_h$, and $p_s^\var$ converges a.e. to $p_s$.
% Finally, \eqref{inCS} implies that $p_h,\;p_s \in L^1([0,T],W^{1,1}(\Omega))$ for all $T>0$.
\medskip 

This concludes the proof of Thm. \ref{thmsimple}.
\bigskip 

{\bf{Proof of Thm. \ref{thmadvanced}}}\\
As in Thm. \ref{thmsimple}, existence (and uniqueness) of strong global solutions to
system (\ref{eqN1bis}) -- (\ref{eqN3bis}), for which $N^\var,\;p^\var_s,\;p^\var_h$ are nonnegative, for a given $\var>0$, is classical (cf. \cite{D_milan}).\medskip
 
Also as in the proof of Thm. \ref{thmsimple}, for each $T>0$, one can find $C_T>0$ such that 
 \begin{equation}\label{nn1}
 \sup_{\var>0} || N^\var ||_{L^{\infty}([0,T]\times \Omega)} \le C_T, 
 \end{equation} 
and as a consequence, there exists $N \in  L^{\infty}([0,T]\times \Omega)$ and a subsequence, still denoted by $N^\var$, such that $N^\var \rightharpoonup N$ in $L^{\infty}([0,T]\times \Omega)$ weak $*$.
\par 
Then, adding (\ref{eqN2bis}) and (\ref{eqN3bis}), we see that (\ref{eqsum1}), (\ref{plu}) still holds, so that using the duality lemma of \cite{CDF}, we end up with
\begin{equation}\label{nn2}
 \sup_{\var>0} || P^\var ||_{L^{2+\delta}([0,T]\times \Omega)} \le C_T , 
  \end{equation}
for some $\delta>0$, and as a consequence, there exist $p_h, p_s \in  L^{2+\delta}([0,T]\times \Omega)$ and subsequences, still denoted by $p_h^\var, p_s^\var$, such that $p_s^\var \rightharpoonup p_s$, $p_h^\var \rightharpoonup p_h$ in $L^{2+\delta}([0,T]\times \Omega)$ weak for some $\delta>0$. 
\medskip

Now observing that the r.h.s. of \eqref{eqN1bis} is bounded in $L^{\infty}([0,T]\times \Omega)$ (this held only in $L^{2+\delta}([0,T]\times \Omega)$ in Thm. \ref{thmsimple}), the maximal regularity estimates for the heat kernel yield the bounds 
\begin{equation}\label{maxreg2} 
\sup_{\var>0} || \pa_t N^\var ||_{L^p([0,T]\times \Omega)} \le C_T,
\qquad
\sup_{\var>0} || \pa_{x_i x_j} N^\var ||_{L^p([0,T]\times \Omega)} \le C_T, 
\end{equation}
for all $p \in [1, +\infty[$, $i,j =1,..,d$, so that the sequence $N^{\var}$ is strongly compact in $L^p([0,T]\times \Omega)$ for all $p \in [1, +\infty[$, and $N^\var$ converges a.e. (up to extraction of a subsequence) towards $N$.
 \medskip

We now compute the derivative of the following nonnegative function:
$$ \frac12  \int \left[ \dfrac{\tilde{\gamma}}{\alpha}(p_h^\var)^2 + N^\var \psi(p_s^\var) \right], $$
with $\psi(x) := \frac{2}{\xi}\left(x - \frac{\ln(1+\xi x)}{\xi}\right)$ (so that $\psi(x) \ge 0$, $\psi'(x) = \frac{2x}{1+\xi x}$, and $\psi''(x) = \frac{2}{(1+\xi x)^2}$). We end up with 
\begin{multline*}
\frac12 \frac{d}{dt}  \int \left[ \dfrac{\tilde{\gamma}}{\alpha}(p_h^\var)^2 + N^\var \psi(p_s^\var) \right]=\\
\int \dfrac{\tilde{\gamma}}{\alpha}p_h^\var \left(d_3\Delta_x p_h^\var-\dfrac{1}{\var}\left( \tilde{\gamma}p_h^\var-\dfrac{\alpha N^\var p_s^\var}{1+\xi p_s^\var}\right)-\mu p_h^\var\right)\\
+\int \dfrac{\psi(p_s^\var)}{2}\left( d_1\Delta_xN^\var+r_0(1-\eta N^\var)N^\var-\dfrac{\alpha N^\var p_s^\var}{1+\xi p_s^\var}\right)\\
+\int N^\var \dfrac{p_s^\var}{1+\xi p_s^\var}\left( d_2 \Delta_x p_s^\var+\dfrac{1}{\var}\left(\tilde{\gamma}p_h^\var-\dfrac{\alpha N^\var p_s^\var}{1+\xi p_s^\var}\right)-\mu p_s^\var+\Gamma p_h^\var\right)
\end{multline*}
%% CALCOLI PER TESI!
%\begin{multline*}
%=\boxed{d_3\dfrac{\tilde{\gamma}}{\alpha}\int{p_h^\var}\Delta_x p_h^\var}
%-\dfrac{1}{\var\alpha}\int{\left(\tilde{\gamma}p_h^\var-\dfrac{\alpha N^\var p_s^\var}{1+\xi p_s^\var}\right)^2}
%-\mu\dfrac{\tilde{\gamma}}{\alpha}\int{(p_h^\var)^2}\\
%+\dfrac{d_1}{2}\int{\psi(p_s^\var)\Delta_x N^\var}
%+\dfrac{r_0}{2}\int{\psi(p_s^\var)(1-\eta N^\var) N^\var}
%-\dfrac{\alpha}{2}\int{\psi(p_s^\var)\dfrac{p_s^\var}{1+\xi p_s^\var}N^\var}\\
%+\boxed{d_2 \int{N^\var\dfrac{p_s^\var}{1+\xi p_s^\var}\Delta_x p_s^\var}}
%-\mu \int{\dfrac{N^\var (p_s^\var)^2}{1+\xi p_s^\var}}
%+\Gamma\int{\dfrac{N^\var p_s^\var p_h^\var}{1+\xi p_s^\var}}
%\end{multline*}
%The boxed terms, integrating by parts and using the boundary conditions, can be rewritten as
%\begin{multline*}
%d_3\dfrac{\tilde{\gamma}}{\alpha}\int{p_h^\var}\Delta_x p_h^\var=-d_3\dfrac{\tilde{\gamma}}{\alpha}\int{|\nabla p_h^\var|^2}\\[-1cm]
%\end{multline*}
%\begin{multline*}
%d_2 \int{N^\var\dfrac{p_s^\var}{1+\xi p_s^\var}\Delta_x p_s^\var}=-d_2\int{\nabla p_s^\var\nabla\left( N^\var\dfrac{p_s^\var}{1+\xi p_s^\var}\right)}\\
%=-d_2\int{\nabla p_s^\var\nabla\left( \dfrac{p_s^\var}{1+\xi p_s^\var}\right)N^\var}-d_2\int{\nabla p_s^\var\nabla N^\var\dfrac{p_s^\var}{1+\xi p_s^\var}}\\
%=-\dfrac{d_2}{2}\int{\nabla p_s^\var\nabla\psi'(p_s^\var)N^\var}-\dfrac{d_2}{2}\int{\nabla \psi(p_s^\var)\nabla N^\var}\\
%=-\dfrac{d_2}{2}\int{|\nabla p_s^\var|^2\nabla\psi''(p_s^\var)N^\var}+\dfrac{d_2}{2}\int{\psi(p_s^\var)\Delta_x N^\var}\\
%\end{multline*}
\begin{multline*}
=-\dfrac{1}{\var\alpha}\underbrace{\int{\left(\tilde{\gamma}p_h^\var-\dfrac{\alpha N^\var p_s^\var}{1+\xi p_s^\var}\right)^2}}_{\boxed{1}} -\underbrace{d_3\dfrac{\tilde{\gamma}}{\alpha}\int{|\nabla_x p_h^\var|^2}}_{\boxed{2}}
-\underbrace{\mu\dfrac{\tilde{\gamma}}{\alpha}\int{(p_h^\var)^2}}_{\boxed{3}}\\
-\underbrace{\dfrac{d_2}{2}\int{N^\var \psi''(p_s^\var)|\nabla_x p_s^\var|^2}}_{\boxed{4}}
+\underbrace{\dfrac{d_2}{2}\int{\psi(p_s^\var)\Delta_x N^\var}}_{\boxed{5}}
-\underbrace{\mu \int{\dfrac{N^\var (p_s^\var)^2}{1+\xi p_s^\var}}}_{\boxed{6}}\\
+\underbrace{\Gamma\int{\dfrac{N^\var p_s^\var p_h^\var}{1+\xi p_s^\var}}}_{\boxed{7}}
+\underbrace{\dfrac{d_1}{2}\int{\psi(p_s^\var)\Delta_x N^\var}}_{\boxed{8}}
+\underbrace{\dfrac{r_0}{2}\int{\psi(p_s^\var)(1-\eta N^\var) N^\var}}_{\boxed{9}}\\
-\underbrace{\dfrac{\alpha}{2}\int{\psi(p_s^\var)\dfrac{p_s^\var}{1+\xi p_s^\var}N^\var}}_{\boxed{10}} .
\end{multline*}
The terms \boxed{1}, \boxed{2}, \boxed{3},  \boxed{4}, \boxed{6} and \boxed{10} are nonpositive. Then remembering that $0\leq \psi(x)\leq \frac{2x}{\xi}$, $p_s^\var$ and $p_h^\var$ are bounded in $L^{2+\delta}$ for some $\delta>0$ (cf. (\ref{nn2})), $N^\var$ is bounded in $L^\infty$
 (cf. (\ref{nn1})), and finally $\Delta_x N^\var$ is bounded in $L^p$ for all $p<+\infty$ (cf. (\ref{maxreg2})), we see that term \boxed{5} and all terms \boxed{7} to \boxed{10}, once integrated on $[0,T]$, are bounded (by some constant $C_T>0$). As a consequence, we end up with the estimates 
\begin{align} 
&\int_0^T\int_{\Omega} \bigg(\tilde{\gamma}p_h^\var - \frac{\alpha N^\var p_s^\var}{1+ \xi p_s^\var} \bigg)^2 \le C_T\var, \label{compj}\\
&\int_0^T\int_{\Omega} |\nabla_x p_h^\var|^2 \le C_T, \label{comp2}\\
&\int_0^T\int_{\Omega} N^\var \psi''(p_s^\var) |\nabla_x p_s^\var|^2 \le C_T.\label{comp}
\end{align}
\medskip

We see that (with $C_T := \frac{\alpha}{\xi} + r_0\eta \sup_{\var>0} ||N^\var||_{L^{\infty} ([0,T]\times \Omega)}$)
$$ (\pa_t - d_1 \Delta_x) N^\var  \ge - C_T N^\var , $$
so that (denoting by $\inf$ the essential infima)
$$ \inf_{\var>0, x\in \Omega} N^\var(t,x) \ge e^{ - C_T} \inf_{\var>0, x\in \Omega} N^\var(0,x). $$

As a consequence, thanks to \eqref{comp},
$$ \int_0^T\int_{\Omega}  \psi''(p_s^\var)  |\nabla_x p_s^\var|^2 \le C_T . $$
Then, Cauchy-Schwarz inequality ensures that 
\begin{equation}
\bigg( \int_0^T\int_{\Omega} |\nabla_x p_s^\var| \bigg)^2 \le  \int_0^T\int_{\Omega}  \psi''(p_s^\var)  |\nabla_x p_s^\var|^2 \int_0^T\int_{\Omega} \frac{(1 + \xi p_s^\var)^2}2 \le C_T. \label{ingps2}
\end{equation}
Using \eqref{eqsum1}, we see that $\pa_t P^\var$ is bounded in $L^{2}([0,T]; H^{-2}( \Omega))$, so that thanks to Aubin's lemma \cite{moussa}, $P^\var$ converges a.e. to $P$. Note then that $ \tilde{\gamma}p_h^\var - \alpha N^\var p_s^\var/(1+\xi p_s^\var) \to 0$ a.e. (up to extraction of a subsequence) thanks to \eqref{compj}, and that $N^\var$ converges a.e. to $N$. Then
$$\tilde{\gamma}P^\var-\left(\tilde{\gamma}p_h^\var-\dfrac{\alpha N^\var p_s^\var}{1+\xi p_s^\var}\right)=\tilde{\gamma}p_s^\var+\dfrac{\alpha N^\var p_s^\var}{1+\xi p_s^\var}$$
converges a.e. towards $\tilde{\gamma}P$. Observing that 
$$\left|\left(\tilde{\gamma}p_s^\var+\dfrac{\alpha N^\var p_s^\var}{1+\xi p_s^\var}\right)-\left(\tilde{\gamma}p_s^\var+\dfrac{\alpha Np_s^\var}{1+\xi p_s^\var}\right)\right|\leq \dfrac{\alpha}{\xi}|N^\var-N|,$$
we see that
$$\tilde{\gamma}p_s^\var+\dfrac{\alpha Np_s^\var}{1+\xi p_s^\var}\to \tilde{\gamma}P.$$
Using the continuity and the strict monotonicity of $y\mapsto \tilde{\gamma}y+\alpha Ny/(1+\xi y)$ for all $N\geq 0$, we see that $p_s^\var$ converges a.e. towards a nonnegative function denoted by $p_s$. Then, $p_h^\var$ also converges a.e. towards a nonnegative function denoted by $p_h$ (because we already know that $P^\var$ converges a.e.). As a consequence, both $p_s^\var$ and  $p_h^\var$ also converge in $L^{2+\delta}([0,T]\times\Omega)$ strong when $\delta>0$ is small enough. Finally, it is clear  that $$p_s+p_h=P, \quad \tilde{\gamma}p_h=\dfrac{\alpha Np_s}{1+\xi p_s},$$
so that \eqref{eqF3bis} holds. We now pass to the limit in equation \eqref{eqN1bis} and \eqref{eqsum1} in the sense of distributions (more precisely, in the sense of very-weak solutions, which include the Neumann boundary conditions and the initial data $N(0,x)=N_{in}(x)$ and 
$(p_s+p_h)(0,x)=p_{s,in}(x)+p_{h,in}(x)$), so that \eqref{eqF1bis} and \eqref{eqF2bis}, or  (\ref{vw1bis}), (\ref{vw2bis}), hold.\medskip

Finally, thanks to estimate (\ref{maxreg2}), we see that $N$ lies in $L^p([0,T],W^{2,p}(\Omega))\cap W^{1,p}(]0,T[,L^p(\Omega))$ for all $p\in]1,+\infty[$, and thanks to estimates \eqref{comp2} and \eqref{ingps2}, we see that $p_s$ and $p_h$ respectively belong to $L^1([0,T],W^{1,1}(\Omega))$ and $L^2([0,T],H^1(\Omega)) \cap L^{\infty}([0,T] \times \Omega)$.

%%%%%%%%%%%%%%%%%%%%%%%%%%%%%%%%%%%%%%%%%%%%%%%%%%%%%%%%%%%%%%%%%%%%%%%%%%%%%%%%%%%%%%%%%%%%%%%%%%%%%%%%%%%%%%%%%%%%%%%%%%%%%%%%%%%%%%
\section{Turing instability analysis}\label{SectionTIa}

\subsection{Limiting system: explicit formulas}\label{limsysder}
In the sequel we systematically assume that $\eta>0$, and use $k :=1/\eta$ in the system \eqref{eqN1bis} -- \eqref{eqN3bis}, so that $k>0$ is the  carrying capacity.  The system becomes
\begin{align}
\dfrac{\partial N^\var}{\partial t} - d_1 \Delta_x N^\var    &= r_0\left(1-\dfrac{N^\var}{k}\right)N^\var- \dfrac{\alpha N^\var p_s^\var}{1+\xi p_s^\var}, \nonumber\\[0.3cm]
\dfrac{\partial p_s^\var}{\partial t} - d_2 \Delta_x p_s^\var&=\dfrac{1}{\varepsilon}\left(-\dfrac{\alpha N^\var p_s^\var}{1+\xi p_s^\var}+\tilde{\gamma} p_h^\var\right)+\Gamma p_h^\var-\mu p_s^\var, \label{qssasystem2}\\[0.3cm]
\dfrac{\partial p_h^\var}{\partial t} - d_3 \Delta_x p_h^\var&=-\dfrac{1}{\varepsilon}\left(-\dfrac{\alpha N^\var p_s^\var}{1+\xi p_s^\var}+\tilde{\gamma} p_h^\var\right)-\mu p_h^\var.\nonumber
\end{align}
We recall that, according to the computation of Section \ref{RigRes}, we know that in the limit when $\var\to 0$, the solution $N^\var,\; p_s^\var,\; p_h^\var$ of this system converges towards $N\geq 0,\; p_s\geq 0,\; p_h\geq 0$ such that 
\begin{equation}
%\tilde{\gamma} p_h=\dfrac{\alpha N p_s}{1+\xi p_s} \qquad \textnormal{so that}\qquad
 \tilde{\gamma}\, p_h= \dfrac{\alpha N p_s}{1+\xi p_s},\label{p_h}\end{equation}
and
$$P:=p_s+p_h=p_s+\dfrac{1}{\tilde{\gamma}}\dfrac{\alpha N p_s}{1+\xi p_s}= \bigg[\dfrac{\tilde{\gamma}\,(1+\xi p_s)+\alpha N}{\tilde{\gamma}(1+\xi p_s)} \bigg]\, p_s.$$
We now wish to write the limiting system 
\begin{align}\begin{split}
\pa_t& N - d_1 \Delta_x N= r_0 \left(1 - \dfrac{N}{k}\right) N - \frac{\alpha N p_s}{1 + \xi p_s}, \\
\pa_t& (p_s + p_h)  -  \Delta_x (d_2 p_s + d_3 p_h) =  \Gamma p_h - \mu (p_h + p_s), \\
\end{split}\label{limsys1}\end{align}
in terms of $N$ and $P$ only. We note that $p_s$ satisfies a second degree equation (when $P$ is given):
%$$P\tilde{\gamma}(1+\xi p_s)=\tilde{\gamma}(1+\xi p_s)p_s+\alpha Np_s$$
$$\tilde{\gamma}\xi p_s^2+(\tilde{\gamma}+\alpha N -\tilde{\gamma}\xi P) p_s - \tilde{\gamma}P=0,$$
so that (considering only the positive root of this equation):
$$p_s=\dfrac{-A+\sqrt{\Delta}}{2\tilde{\gamma}\xi}=\dfrac{2\tilde{\gamma}P}{A+\sqrt{\Delta}},$$
where we have denoted
%$$\Delta=(\tilde{\gamma}+\alpha N -\tilde{\gamma}\xi P)^2+4\tilde{\gamma}^2\xi P>0$$
\begin{equation}
A:=\tilde{\gamma}+\alpha N -\tilde{\gamma}\xi P,\qquad \Delta=A^2+4\tilde{\gamma}^2\xi P.\label{ADelta}
\end{equation}
Note that $\Delta>0$ since $P>0$.
Denoting by 
\begin{equation}
B:=\tilde{\gamma}+\alpha N +\tilde{\gamma}\xi P,   \label{B}
\end{equation}
from \eqref{p_h} we also obtain
$$p_h=\dfrac{2\alpha NP}{B+\sqrt{\Delta}},$$
where $\Delta$ can be computed in terms of $B$:
\begin{align*}
\Delta&=A^2+4\tilde{\gamma}^2\xi P\\
%=(\tilde{\gamma}+\alpha N -\tilde{\gamma}\xi P)^2+4\tilde{\gamma}^2\xiP\\
%      &=\tilde{\gamma}^2+(\alpha N)^2+(\tilde{\gamma}\xiP)^2+2\tilde{\gamma}\alpha N-2\tilde{\gamma}^2\xi P-2\alpha N \tilde{\gamma}\xi P+4\tilde{\gamma}^2\xi P\\
	%&=\tilde{\gamma}^2+(\alpha N)^2+(\tilde{\gamma}\xi P)^2+2\tilde{\gamma}\alpha N+2\tilde{\gamma}^2\xi P-2\alpha N \tilde{\gamma}\xi P=\\
%			&=\tilde{\gamma}^2+(\alpha N)^2+(\tilde{\gamma}\xi P)^2+2\tilde{\gamma}\alpha N+2\tilde{\gamma}^2\xi P+2\alpha N \tilde{\gamma}\xi P-4\alpha N \tilde{\gamma}\xi P\\
	%		&=(\tilde{\gamma}+\alpha N +\tilde{\gamma}\xi P)^2-4\alpha N \tilde{\gamma}\xi P
	&=B^2-4\alpha N \tilde{\gamma}\xi P.
\end{align*}
%Adding the second and the third equation of \eqref{qssasystem2}
%$$\dfrac{\partial P}{\partial t} - \Delta_x (d_2 p_s+d_3 p_h)=\Gamma p_h-\mu P$$
%$$\dfrac{\partial P}{\partial t} - \Delta_x \left(d_2 \dfrac{2\tilde{\gamma}P}{A+\sqrt{\Delta}}+d_3 \dfrac{2\alpha NP}{B+\sqrt{\Delta}}\right)=\Gamma\dfrac{2\alpha  NP}{B+\sqrt{\Delta}}-\mu P$$
Then the limiting system can be written with $N,\; P$ as unknowns in the following way:
\begin{align}\begin{split}
\dfrac{\partial N}{\partial t}& - d_1 \Delta_x N = r_0\left(1-\dfrac{N}{k}\right)N- \tilde{\gamma}\dfrac{2\alpha  NP}{B+\sqrt{\Delta}},\\[0.3cm]
\dfrac{\partial P}{\partial t}& - \Delta_x \left(d_2 \dfrac{2\tilde{\gamma}P}{A+\sqrt{\Delta}}+d_3 \dfrac{2\alpha NP}{B+\sqrt{\Delta}}\right)=\Gamma\dfrac{2\alpha  NP}{B+\sqrt{\Delta}}-\mu P,
\end{split}\label{qssasystem3}\end{align}
where  $A,\; B$ and $\Delta$ are defined in \eqref{ADelta} and \eqref{B}.

\subsection{Adimensionalization}

In order to simplify the notations and to keep only meaningful parameters, we now propose an adimensionalization procedure, using the new variables $T,\; n,\; p$ instead of $t,\; N,\; P$ in the following way: 
$$t=\Theta T,\; N=\Sigma n,\; P=\Pi p.$$
%$$\begin{cases}
%\dfrac{\Sigma}{\Theta}\dfrac{\partial n}{\partial T} - d_1\Sigma \Delta_x n = r_0\left(1-\dfrac{\Sigma n}{k}\right)\Sigma n- \tilde{\gamma}\dfrac{2\alpha \Sigma \Pi np}{B+\sqrt{\Delta}}\\[0.3cm]
%\dfrac{\Pi}{\Theta}\dfrac{\partial p}{\partial T} - \Delta_x \left(d_2 \dfrac{2\tilde{\gamma}\Pi p}{A+\sqrt{\Delta}}+d_3 \dfrac{2\alpha \Sigma \Pi np}{B+\sqrt{\Delta}}\right)=\Gamma\dfrac{2\alpha\Sigma\Pi np}{B+\sqrt{\Delta}}-\mu \Pi p
%\end{cases}$$
%or simplifying
After simplifications, the system \eqref{qssasystem3} becomes
\begin{align*}
\dfrac{\partial n}{\partial T}& - d_1\Theta \Delta_x n = r_0\Theta\left(1-\dfrac{n}{k/\Sigma}\right) n- \dfrac{2\tilde{\gamma}\alpha\Pi\Theta np}{B+\sqrt{\Delta}},\\[0.3cm]
\dfrac{\partial p}{\partial T}& - \Delta_x \left(d_2\Theta \dfrac{2\tilde{\gamma}p}{A+\sqrt{\Delta}}+d_3\Theta \dfrac{2\alpha\Sigma np}{B+\sqrt{\Delta}}\right)=\dfrac{2\Gamma\alpha\Sigma\Theta np}{B+\sqrt{\Delta}}-\mu\Theta p,
\end{align*}
where
\begin{align*}
A+\sqrt{\Delta}&=\tilde{\gamma}+\alpha \Sigma n -\tilde{\gamma}\xi \Pi p+\sqrt{(\tilde{\gamma}+\alpha \Sigma n +\tilde{\gamma}\xi \Pi p)^2-4\tilde{\gamma}\alpha\xi\Sigma\Pi np},\\
B+\sqrt{\Delta}&=\tilde{\gamma}+\alpha \Sigma n +\tilde{\gamma}\xi \Pi p+\sqrt{(\tilde{\gamma}+\alpha \Sigma n +\tilde{\gamma}\xi \Pi p)^2-4\tilde{\gamma}\alpha\xi\Sigma\Pi np}.
\end{align*}
Choosing $\Theta,\; \Sigma,\; \Pi$ in such a way that $2\alpha\Pi\Theta=1,\; \alpha\Sigma=1,\; \tilde{\gamma}\xi\Pi=1$, we end up with the system
\begin{align*}
&\dfrac{\partial n}{\partial T} - d_1\Theta \Delta_x n = r_0\Theta\left(1-\dfrac{n}{k/\Sigma}\right) n- \dfrac{\tilde{\gamma} np}{B+\sqrt{\Delta}},\\[0.3cm]
&\dfrac{\partial p}{\partial T} - \Delta_x \left(d_2\Theta \dfrac{2\tilde{\gamma}p}{A+\sqrt{\Delta}}+d_3\Theta \dfrac{2 np}{B+\sqrt{\Delta}}\right)=\dfrac{\Gamma\tilde{\gamma}\Sigma\xi np}{B+\sqrt{\Delta}}-\mu\Theta p,
\end{align*}
where now
\begin{align*}
A+\sqrt{\Delta}&=\tilde{\gamma}+n-p+\sqrt{(\tilde{\gamma}+n+p)^2-4np},\\
B+\sqrt{\Delta}&=\tilde{\gamma}+n+p+\sqrt{(\tilde{\gamma}+n+p)^2-4np}.
\end{align*}
We set $D_1:=d_1\Theta,\; D_2:=d_2\Theta,\; D_3:=d_3\Theta,\; r=r_0\Theta,\; \nu:=k/\Sigma$.\\
%{\dafare Problema per $\Gamma$ nell'adimensionalizzazione???}\\
Furthermore, we denote again $n$ by $N$, $p$ by $P$, $\tilde{\gamma}$ by $\gamma$, and redefine $\Gamma:=\Gamma\tilde{\gamma}\Sigma\xi,\;\mu:=\mu\Theta$. We end up with
\begin{align}\begin{split}
\dfrac{\partial N}{\partial T}& - D_1\Delta_x N = r\left(1-\dfrac{N}{\nu}\right)N-\dfrac{\gamma NP}{B+\sqrt{\Delta}},\\[0.3cm]
\dfrac{\partial P}{\partial T}& - \Delta_x \left(D_2\dfrac{2\gamma P}{A+\sqrt{\Delta}}+D_3\dfrac{2 NP}{B+\sqrt{\Delta}}\right)=\dfrac{\Gamma NP}{B+\sqrt{\Delta}}-\mu P,
\end{split}\label{limsys}\end{align}
where now
\begin{equation}
A=\gamma+N-P,\quad B=\gamma+N+P,\quad \Delta=(\gamma+N+P)^2-4NP.\label{ABDelta}
\end{equation}

Rationalizing the denominators, we can obtain an equivalent expression, which is useful for the stability analysis of equilibrium states:
%\begin{align}\begin{split}
%\dfrac{\partial N}{\partial T}& - D_1\Delta_x N = r\left(1-\dfrac{N}{\nu}\right)N-\dfrac{\gamma}{4}\left(\gamma+N+P-\sqrt{(\gamma+N+P)^2-4NP}\right),\nonumber\\[0.3cm]
%\dfrac{\partial P}{\partial T}& - \Delta_x \left(\dfrac{D_2}{2}{(\sqrt{\Delta}-A)}+\dfrac{D_3}{2}{(B-\sqrt{\Delta})}\right)=\\
%&\dfrac{\Gamma}{4}\left(\gamma+N+P-\sqrt{(\gamma+N+P)^2-4NP}\right)-\mu P,
%\end{split}\label{limsysratio}\end{align}
\begin{align}\begin{split}
\dfrac{\partial N}{\partial T}& - D_1\Delta_x N = r\left(1-\dfrac{N}{\nu}\right)N-\dfrac{\gamma}{4}\left(B-\sqrt{B^2-4NP}\right),\\[0.3cm]
\dfrac{\partial P}{\partial T}& - \Delta_x \left(\dfrac{D_2}{2}{(\sqrt{\Delta}-A)}+\dfrac{D_3}{2}{(B-\sqrt{\Delta})}\right)=
\dfrac{\Gamma}{4}\left(B-\sqrt{B^2-4NP}\right)-\mu P,
\end{split}\label{limsysratio}\end{align}
where $A,\; B$ and $\Delta$ are still defined by \eqref{ABDelta}.

The limiting system presents a cross-diffusion term in the predator equation (the diffusion rate depends on the prey biomass), and a trophic function close to the Beddington-DeAngelis one. 
\medskip

\subsection{Homogeneous equilibrium states}

In this subsection, we look for the equilibrium states of the ODEs system corresponding to the reaction part of the whole system \eqref{limsys}, or equivalently \eqref{limsysratio}:
\begin{align}\begin{split}
\dot{N}&=r\left(1-\dfrac{N}{\nu}\right)N-\dfrac{\gamma NP}{B+\sqrt{\Delta}},\\[0.5cm]
\dot{P}&=\dfrac{\Gamma NP}{B+\sqrt{\Delta}}-\mu P,
\end{split}\label{limsyseq}\end{align}
where $A,\; B$ and $\Delta$ are defined in \eqref{ABDelta}.\medskip

From the first equation, if $P=0$ we obtain $N=0$ or $N=\nu$, corresponding to the total extinction $E_0(0,0)$ and the non-coexistence $E_1(\nu,0)$ equilibria. \medskip 

Otherwise, we look for a coexistence equilibrium $E_*(N_*,P_*)$ (that means $N_*,\;P_*>0$). From the second equation, we get the identity
\begin{equation}
\dfrac{\Gamma N_*}{\gamma+N_*+P_*+\sqrt{(\gamma+N_*+P_*)^2-4N_*P_*}}-\mu=0, 
\label{eq per eq} \end{equation}
from which, rationalizing the denominator, we can obtain 
\begin{equation}
\sqrt{(\gamma+N_*+P_*)^2-4N_*P_*}=\gamma+N_*+P_*-\dfrac{4\mu}{\Gamma}P_*.
\label{sqrt*}\end{equation}
Rewriting \eqref{eq per eq} as 
%, we can also derive an expression of $P_*$ in terms of $N_*$. In fact, we have 
$$\Gamma N_*-\mu(\gamma+N_*+P_*)=\mu \sqrt{(\gamma+N_*+P_*)^2-4N_*P_*},$$
we see that the searched equilibrium can exist only if
\begin{equation}
\Gamma N_*-\mu(\gamma+N_*+P_*)>0.
\label{cond1}\end{equation}
Taking the square of both terms, we end up with
\begin{equation}
P_*=\dfrac{\Gamma}{2\mu}\dfrac{(\Gamma-2\mu)N_*-2\gamma\mu}{(\Gamma-2\mu)}=\dfrac{\Gamma}{2\mu}N_*-\dfrac{\Gamma\gamma}{\Gamma-2\mu}.
\label{P*N*}\end{equation}
Substituting \eqref{P*N*} in \eqref{cond1}, we see that \eqref{cond1} is equivalent to 
\begin{equation} 
\dfrac{\Gamma-2\mu}{2}N_*+\gamma\mu\left(\dfrac{2\mu}{\Gamma-2\mu}\right)>0. 
\label{cond1mod}\end{equation}
Since we are looking for equilibria with $N_*>0$. we see that \eqref{cond1} or \eqref{cond1mod} can also be rewritten as
\begin{equation} 
\Gamma-2\mu>0. 
\label{cond1mod2}\end{equation}
Substituting the expression \eqref{sqrt*} in the equation $\dot{N}=0$ from the first equation of \eqref{limsyseq} written as
$$\dot{N}=r\left(1-\dfrac{N}{\nu}\right)N-\dfrac{\gamma}{4}{(B-\sqrt{\Delta})},$$
we obtain 
\begin{equation}
r\left(1-\dfrac{N}{\nu}\right)N-\dfrac{\gamma\mu}{\Gamma}P_*=0,
\label{eqforN*}\end{equation}
from which we have another expression of $P_*$ in terms of $N_*$:
\begin{equation}
P_*=\dfrac{\Gamma}{\gamma\mu} r\left(1-\dfrac{N_*}{\nu}\right)N_*.
\label{P*N*2}\end{equation}
Substituting the expression \eqref{P*N*} in \eqref{eqforN*}, we obtain a second order equation in the unknown $N_*$:
\begin{equation}
\dfrac{r}{\nu}N_*^2-\left(r-\dfrac{\gamma}{2}\right)N_*-\dfrac{\mu\gamma^2}{\Gamma-2\mu}=0.
\label{2ndorderN*}\end{equation}
We see that, thanks to \eqref{cond1mod2},
\begin{equation}
\Delta_N:=\left(r-\dfrac{\gamma}{2}\right)^2+4\dfrac{r}{\nu}\dfrac{\mu\gamma^2}{\Gamma-2\mu}>0,
\label{cond2}\end{equation}
so that equation \eqref{2ndorderN*} has one and only one strictly positive solution,
 given by
\begin{equation} 
N_{*}=\dfrac{\nu}{2r}\left(r-\dfrac{\gamma}{2}+\sqrt{\Delta_N}\right). 
\label{N*pm} \end{equation}
Then, the condition $P_*>0$ is equivalent to
\begin{equation}
\dfrac{\Gamma}{2\mu}N_*-\dfrac{\Gamma\gamma}{\Gamma-2\mu}>0 \;\Leftrightarrow \; 0<N_*<\nu,
\label{cond4}\end{equation}
depending on the chosen expression for $P_*$.
This condition can be rewritten as
\begin{equation}%\tag{COND $E_*$} 
\Gamma-2\mu>\dfrac{2\gamma\mu}{\nu},
\label{COND_Estar}\end{equation}
by substituting \eqref{N*pm} in the last term of \eqref{cond4}. Note that this last necessary condition for the existence of the coexistence equilibrium $E_*$ implies condition \eqref{cond1mod2}. We now briefly explain why condition (\ref{COND_Estar}) is in fact both necessary and sufficient for the existence of $E_*$. Indeed, (\ref{COND_Estar}) can be rewritten as 
$$\sqrt{\Delta_N}<r+\dfrac \gamma 2,$$
so that $N_*$ (computed from formula \eqref{N*pm} and \eqref{cond2}) is such that $0<N_*<\nu$. Remembering that this last condition is equivalent to $P_*>0$ (when $P_*$ is given by \eqref{P*N*2} for example), we see that both $N_*$ and $P_*$ defined in this way are strictly positive. One can easily check that they satisfy $\dot{N}=0,\; \dot{P}=0$ in \eqref{limsyseq}.
\bigskip

%%%%%%%%%%%%%%%
We now study the stability properties of these equilibrium states. We denote as $J_{ij},\; i,j=1,2$, the elements of the Jacobian matrix of the system \eqref{limsyseq}:
\begin{align*}
J_{11}&=\dfrac{\partial}{\partial N} \dot{N}=r-\dfrac{2r}{\nu}N-\dfrac{\gamma}{4}\left(1-\dfrac{\gamma+N-P}{\sqrt{(\gamma+N+P)^2-4NP}}\right),\\
J_{12}&=\dfrac{\partial}{\partial P} \dot{N}=-\dfrac{\gamma}{4}\left(1-\dfrac{\gamma-N+P}{\sqrt{(\gamma+N+P)^2-4NP}}\right),\\
J_{21}&=\dfrac{\partial}{\partial N} \dot{P}=\dfrac{\Gamma}{4}\left(1-\dfrac{\gamma+N-P}{\sqrt{(\gamma+N+P)^2-4NP}}\right),\\
J_{22}&=\dfrac{\partial}{\partial P} \dot{P}=\dfrac{\Gamma}{4}\left(1-\dfrac{\gamma-N+P}{\sqrt{(\gamma+N+P)^2-4NP}}\right)-\mu.
\end{align*}
Evaluating the Jacobian matrix in the equilibrium states, we obtain:
$$J(E_0)=\begin{bmatrix} r & 0 \\ 0 &-\mu \end{bmatrix},\qquad J(E_1)=\begin{bmatrix} -r & * \\ 0 &J_{22}(E_1) \end{bmatrix}$$
(where $*$ in the matrix means ``some term that we do not make explicit''), and
$$J_{22}(E_1)=\dfrac{\Gamma\nu}{2(\gamma+\nu)}-\mu>0 \qquad \Leftrightarrow \qquad \Gamma-2\mu>\dfrac{2\gamma\mu}{\nu}.$$
Then the equilibrium $E_0$ is unstable (it is a saddle point); the equilibrium $E_1$ is locally asymptotically stable (it is a node) when $E_*$ does not exist, and unstable (it is a saddle point) otherwise.
\medskip

The study of the stability of $E_*$ is more intricate. First, we compute the quantities
\begin{equation}
Q_*:=\gamma+N_*+P_*-\dfrac{4\mu}{\Gamma}P_*=\dfrac{\Gamma-2\mu}{2\mu}N_*+\dfrac{2\mu\gamma}{\Gamma-2\mu}>0
\label{Qstar}\end{equation}
and, thank to \eqref{P*N*} and \eqref{Qstar},
$$\dfrac{P_*}{Q_*}=\dfrac{\Gamma}{\Gamma-2\mu}\left(1-\dfrac{2\mu\gamma\Gamma}{(\Gamma-2\mu)^2N_*+4\mu^2\gamma}\right).$$
Thanks to \eqref{sqrt*} and \eqref{P*N*}, we obtain (remembering the definition of $\Delta_N$ in \eqref{cond2}) explicit expressions for the elements of the Jacobian matrix evaluated at the equilibrium $E_*$, denoted by $\bm{J}^* := J(E_*)$:
\begin{align*}
J^*_{11}:=J_{11}(E_*)&=r\left(1-\dfrac{2}{\nu}N_*\right)-\dfrac{\gamma(\Gamma-2\mu)}{2\Gamma}\dfrac{P_*}{Q_*}\\
           &=-\dfrac{r}{Q_*}\dfrac{2\mu\gamma}{\Gamma-2\mu}\left[\dfrac{(\Gamma-2\mu)^2}{4\nu\mu^2\gamma}N_*^2+\dfrac{2}{\nu}N_*-1 \right]\\
           &=\dfrac{\mu\gamma^2\Gamma}{(\Gamma-2\mu)^2N_*+4\mu^2\gamma}-\sqrt{\Delta_N},\\[0.3cm]
J^*_{12}:=J_{12}(E_*)&=-\dfrac{\gamma}{4}\left(\dfrac{\sqrt{(\gamma+N_*+P_*)^2-4N_*P_*}-(\gamma+N_*+P_*)+2N_*}{\sqrt{(\gamma+N_*+P_*)^2-4N_*P_*}}\right)\\
&=-\dfrac{\gamma}{2\Gamma Q_*}(-2\mu P_*+\Gamma N_*)=-\dfrac{\gamma^2\mu}{Q_*(\Gamma-2\mu)}<0,\\[0.3cm]
J^*_{21}:=J_{21}(E_*)&=\dfrac{\Gamma}{4}\left(\dfrac{\sqrt{(\gamma+N_*+P_*)^2-4N_*P_*}-(\gamma+N_*+P_*)+2P_*}{\sqrt{(\gamma+N_*+P_*)^2-4N_*P_*}}\right)\\
&=\dfrac{1}{2Q_*}(\Gamma-2\mu)P_*>0,\\[0.3cm]
J^*_{22}:=J_{22}(E_*)&=-\mu+\dfrac{\Gamma}{4}\left(\dfrac{\sqrt{(\gamma+N_*+P_*)^2-4N_*P_*}-(\gamma+N_*+P_*)+2N_*}{\sqrt{(\gamma+N_*+P_*)^2-4N_*P_*}}\right)\\
&=-\mu+\dfrac{\Gamma}{2Q_*}\left(N_*-\dfrac{2\mu}{\Gamma}P_*\right)\\
&=\dfrac{1}{2Q}\left(-2\mu\left(\gamma+N_*+P_*-\dfrac{4\mu}{\Gamma}P_*\right)+\Gamma N_*-2\mu P_*\right)\\
          &=\dfrac{1}{2Q_*}\left((-2\mu\gamma +(\Gamma-2\mu) N_*) - 4\mu P_*+\dfrac{8\mu^2}{\Gamma}P_*  \right)\\
					&=\dfrac{1}{2Q_*}\left(\dfrac{2\mu}{\Gamma}(\Gamma-2\mu)P_*- 4\mu P_*+\dfrac{8\mu^2}{\Gamma}P_*  \right)=-\dfrac{\mu}{\Gamma Q_*}(\Gamma-2\mu)P_*<0.
\end{align*}

Note that the sign of $J^*_{11}$ is not prescribed, while we are able to determine the sign of all the others elements ($J^*_{12}<0,\; J^*_{21}>0,\; J^*_{22}>0$). However, we are able to prove that 
$$\det \bm{J}^*=J^*_{11}J^*_{22}-J^*_{12}J^*_{21}>0,$$
whatever is the sign of $J^*_{11}$. In fact, substituting in the expression of $\det \bm{J}^*$ the formulas of $J^*_{ij},\; i,j=1,2$, we get
%\begin{multline*}
$$
\det \bm{J}^*
%=\dfrac{2\mu r^2}{Q_*^2}\left(\dfrac{(\Gamma-2\mu)^2}{4\nu\mu^2\gamma}N_*^2+\dfrac{2}{\nu}N_*-1\right)\left(1-\dfrac{N_*}{\nu}\right) N_*+\dfrac{\Gamma r}{2Q_*^2}\left(1-\dfrac{N_*}{\nu}\right) N_*\\
=\dfrac{r}{Q_*^2}\left(1-\dfrac{N_*}{\nu}\right) N_*\left[r\left(\dfrac{(\Gamma-2\mu)^2}{2\nu\mu\gamma}N_*^2+\dfrac{4\mu}{\nu}N_*-2\mu \right)+\dfrac{\Gamma \gamma}{2}\right],
$$
%\end{multline*}
and substituting \eqref{N*pm} in the linear term inside the brackets, we obtain
$$
\det \bm{J}^*=\dfrac{r}{Q_*^2}\left(1-\dfrac{N_*}{\nu}\right) N_*\left[r\dfrac{(\Gamma-2\mu)^2}{2\nu\mu\gamma}N_*^2+2\mu\sqrt{\Delta_N} +\dfrac{\gamma}{2}(\Gamma-2\mu)\right]>0.$$
On the contrary, the sign of the trace of the Jacobin matrix evaluated at $E_*$ is not prescribed. Indeed, when $J^*_{11}<0$, we have
$$\tr \bm{J}^* =J^*_{11}+J^*_{22}<0 , $$
and then $E_*$ is locally asymptotically stable. However, when $J^*_{11}>0$, the trace can be nonpositive or nonnegative. Numerical evidences show that both cases can hold for different values of the parameters of the model.

\subsection{Turing Instability when linear diffusions are added to the system of ODEs} \label{sub34}

We consider the system of reaction-diffusion defined by (for given $D_1,D_P>0$)
\begin{align}\begin{split}
\dfrac{\partial N}{\partial t}& - D_1 \Delta_x N= r\left(1-\dfrac{N}{\nu}\right)N-\dfrac{\gamma NP}{B+\sqrt{\Delta}},\\[0.3cm]
\dfrac{\partial P}{\partial t}& - D_P \Delta_x P=\dfrac{\Gamma NP}{B+\sqrt{\Delta}}-\mu P,
\end{split}\label{BdAlinear}\end{align}
where $A,\; B$ and $\Delta$ are defined in \eqref{ABDelta}, and sets of parameter values such that $\tr J(E_*)<0$ (Note that such parameters indeed exist).\\
For any $\lambda_k\ge 0$ eigenvalue of the Laplacian on $\Omega$ (with Neumann boundary conditions), where $k\in \mathbb{N}$, the characteristic matrix evaluated at the equilibrium $E_*$
$$\bm{M}=\begin{bmatrix}J^*_{11}-D_1 \lambda_k & J^*_{12}\\ J^*_{21} & J^*_{22}-D_P \lambda_k\end{bmatrix}$$
has a strictly negative trace. In fact, 
$$ Tr \bm{M}=J^*_{11}-D_1 \lambda_k+J^*_{22}-D_P \lambda_k=\tr \bm{J}-D_1 \lambda_k-D_2 \lambda_k<0.$$ 
Its determinant is
$$\det \bm{M}=\det \bm{J}-(D_PJ^*_{11}+D_1J^*_{22})\lambda_k+D_1D_P\lambda_k^2,$$
so that a necessary condition for the Turing instability to appear is 
$$D_1J^*_{22}+D_PJ^*_{11}>0.$$
If $J^*_{11}<0$, remembering that $J_{22}^*<0$, we see that $D_PJ^*_{11}+D_1J^*_{22}<0$, so that no  Turing instability can appear.

On the opposite, if $J^*_{11}>0$, for any given $k\neq 0$ (so that $\lambda_k>0$), we can select $D_P$ sufficiently large for getting  $\det \bm{J}-D_PJ^*_{11}\lambda_k<0$. Then, when $D_1>0$ is small enough, $\det \bm{M}<0$ and the Turing instability appears.

\subsection{Turing Instability with cross diffusion}

 We now consider system \eqref{limsys} or the equivalent form \eqref{limsysratio}, that is 
\begin{align*}\begin{split}
\dfrac{\partial N}{\partial T}& - D_1\Delta_x N = r\left(1-\dfrac{N}{\nu}\right)N-\dfrac{\gamma}{4}\left(B-\sqrt{B^2-4NP}\right),\\[0.3cm]
\dfrac{\partial P}{\partial T}& - \Delta_x \left(\dfrac{D_2}{2}{(\sqrt{\Delta}-A)}+\dfrac{D_3}{2}{(B-\sqrt{\Delta})}\right)=
\dfrac{\Gamma}{4}\left(B-\sqrt{B^2-4NP}\right)-\mu P,
\end{split}\end{align*}
where $A,\; B$ and $\Delta$ are defined by \eqref{ABDelta}.
The characteristic matrix takes the form
$$\bm{M} :=\begin{bmatrix} J^*_{11}-J^*_{\Delta 11}\lambda_k & J^*_{12}-J^*_{\Delta 12}\lambda_k\\[0.2cm] J^*_{21}-J^*_{\Delta 21}\lambda_k & J^*_{22}-J^*_{\Delta 22}\lambda_k\end{bmatrix}, $$
where the $\lambda_k$ are defined as in Subsection \ref{sub34}, and the terms $J_{\Delta ij},\;i,j=1,2$ are obtained by linearizing the diffusion terms around $E_*$. Their explicit form is given by the following formulas:
\begin{align}
J^*_{\Delta 11}&=D_1>0, \nonumber\\[0.1cm]
J^*_{\Delta 12}&=0,\nonumber\\[0.1cm]
J^*_{\Delta 21}&=\dfrac{D_2-D_3}{2}\left(\dfrac{(\gamma+N_*-P_*)-\sqrt{\Delta_*}}{\sqrt{\Delta_*}} \right)=-\dfrac{D_2-D_3}{Q_*}\dfrac{\Gamma-2\mu}{\Gamma}P_*,\label{JDelta21}\\[0.1cm]
J^*_{\Delta 22}&=\dfrac{D_2}{2}\left(\dfrac{(\gamma-N_*+P_*)+\sqrt{\Delta_*}}{\sqrt{\Delta_*}} \right)+\dfrac{D_3}{2}\left(\dfrac{\sqrt{\Delta_*}-(\gamma-N_*+P_*)}{\sqrt{\Delta_*}} \right)\nonumber\\
&=\dfrac{D_2}{Q_*}\left(\gamma+P_*\dfrac{(\Gamma-2\mu)}{\Gamma}\right)+\dfrac{D_3}{Q_*}\dfrac{2\mu\gamma}{\Gamma-2\mu}>0.\label{JDelta22}
\end{align}
We notice that only the sign of $J^2_{\Delta 21}$ depends on $D_2, D_3$. Due to the biological meaning of these parameters, we systematically assume that $D_2>D_3$. Indeed, searching predators are expected to diffuse more quickly  than handling predators. With this assumption, $J^*_{\Delta 21}<0$.\medskip

We still consider sets of parameter values such that $\tr J(E_*)<0$. Then the characteristic matrix $\bm{M}$ has strictly negative trace, because 
$$\tr \bm{M}=J^*_{11}-D_1\lambda_k+J^*_{22}-J^*_{\Delta 22}\lambda_k=\underbrace{\tr \bm{J}^*}_{-}-D_1 \lambda_k-\underbrace{J^*_{\Delta 22}}_{+} \lambda_k<0.$$  
Its determinant is
$$\det \bm{M}=\underbrace{\det \bm{J}^*}_{+}-(J^*_{11}J^*_{\Delta 22}+J^*_{22}J^*_{\Delta 11}-J^*_{12}J^*_{\Delta 21})\lambda_k+\underbrace{(D_1J^*_{\Delta 22})}_{+}\lambda_k^2 , $$
so that a necessary condition for the Turing instability to appear is 
\begin{equation}
J^*_{11}J^*_{\Delta 22}+J^*_{22}J^*_{\Delta 11}-J^*_{12}J^*_{\Delta 21}>0.
\label{NcTicross}\end{equation}
If $J^*_{11}<0$, we have
$$\underbrace{J^*_{11}}_{-}\underbrace{J^*_{\Delta 22}}_{+}+\underbrace{J^*_{22}}_{-}\underbrace{J^*_{\Delta 11}}_{+}-\underbrace{J^*_{12}}_{-}\underbrace{J^*_{\Delta 21}}_{-},$$
so that condition \eqref{NcTicross} does not hold and, as in the case of linear diffusion, no  Turing instability can appear. On the opposite, if $J^*_{11}>0$, we have
$$\underbrace{J^*_{11}}_{+}\underbrace{J^*_{\Delta 22}}_{+}+\underbrace{J^*_{22}}_{-}\underbrace{J^*_{\Delta 11}}_{+}-\underbrace{J^*_{12}}_{-}\underbrace{J^*_{\Delta 21}}_{-} .$$
%so that the necessary condition for the Turing instability could be satisfied for sufficiently large $D_2$ and $D_3$.
Then, for any $k \neq 0$ (so that $\lambda_k > 0$, we can select $D_2$ large enough and $D_3\sim D_2$ so that $\det \bm{J}-(J^*_{11}J^*_{\Delta 22} -J^*_{12}J^*_{\Delta 21})\lambda_k<0$. Then, when $D_1>0$ is small enough, $\det \bm{M}<0$ and the Turing instability appears.

\subsection{Turing instability regions: linear versus cross diffusion}

We recall that the derivation of equations \eqref{limsys} produces a cross diffusion term in the predator equation, whereas the prey diffusion rate is still a constant. We recall, for reader's convenience, the model equations 
\begin{align}\begin{split}
\dfrac{\partial N}{\partial T}& - D_1\Delta_x N = r\left(1-\dfrac{N}{\nu}\right)N-\dfrac{\gamma NP}{B+\sqrt{\Delta}},\\[0.3cm]
\dfrac{\partial P}{\partial T}& - \Delta_x \left(f(P,N)P\right)=\dfrac{\Gamma NP}{B+\sqrt{\Delta}}-\mu P,
\end{split}\label{sysbis}\end{align}
where
$$A=\gamma+N-P,\quad B=\gamma+N+P,\quad \Delta=(\gamma+N+P)^2-4NP . $$
The term $\Delta_x\left( f(P,N) P\right)$ is the cross diffusion term, and
\begin{equation}
f(P,N):=D_2 \dfrac{2\gamma}{A+\sqrt{\Delta}}+D_3 \dfrac{2N}{B+\sqrt{\Delta}} . \label{FCrossDiffusionTerm}
\end{equation}

We want to compare three natural strategies to model the diffusion in the predator-prey interactions. 
\begin{enumerate}
\item First, we take the reaction part of \eqref{sysbis} and we add a diffusion term with $D_2$ as constant rate for predators. This means that we exactly take the diffusion coefficient of searching predators for all predators appearing in the limiting model.
\item Secondly, we take the reaction part of \eqref{sysbis} and we add a diffusion term with $D_P=f(P_*,N_*)$, where $f$ is defined in \eqref{FCrossDiffusionTerm}, as constant rate of diffusion in the equation for predators. This means that we now take into account the difference among searching and handling predators, since both diffusion rates $D_2$ and $D_3$ are present in equation \eqref{FCrossDiffusionTerm}. 
\item Finally, we consider the cross diffusion model \eqref{sysbis}, coming out of the derivation of the model by singular perturbation as explained in Section \ref{RigRes} and Subsection \ref{limsysder}.
\end{enumerate}
\bigskip

Note first that, thanks to \eqref{P*N*}, it is possible to obtain a simple expression for $D_P$. Indeed,
\begin{align}
D_P&=f(P_*,N_*)=D_2 \dfrac{2\gamma}{A_*+\sqrt{\Delta_*}}+D_3 \dfrac{2N_*}{B_*+\sqrt{\Delta_*}}\nonumber\\
   &=D_2 \dfrac{\gamma}{\gamma+N_*-\dfrac{2\mu}{\Gamma}P_*}+D_3 \dfrac{N_*}{\gamma+N_*+P_*-\dfrac{2\mu}{\Gamma}P_*}\nonumber\\
	 &=D_2 \dfrac{\gamma}{\gamma+\dfrac{2\mu\gamma}{\Gamma-2\mu}}+D_3 \dfrac{N_*}{\gamma+\dfrac{\Gamma}{2\mu}N_*-\dfrac{\Gamma\gamma}{\Gamma-2\mu}+\dfrac{2\mu\gamma}{\Gamma-2\mu}}\nonumber\\
	&=D_2 \left( 1-\dfrac{2\mu}{\Gamma}\right)+D_3\dfrac{2\mu}{\Gamma}.\label{defDP}
\end{align}
We notice therefore that $D_P$ is a convex combination of the diffusion coefficients $D_2,\;D_3$ of searching and handling predators. Furthermore, assuming that $D_3<D_2$ (remember that from the modeling point of view, handling predators have a lower diffusion rate than searching predators), we get the estimate
$$D_P=D_2 \left( 1-\dfrac{2\mu}{\Gamma}\right)+D_3\dfrac{2\mu}{\Gamma}=D_2-\dfrac{D_2-D_3}{\Gamma}2\mu<D_2.$$
The characteristic matrices of the cases that we consider are finally given by
\begin{enumerate}
	\item Linear diffusion with rate $D_2$: 
	\begin{equation}
	\bm{M}_{L2}=\begin{bmatrix} J^*_{11}-D_1\lambda_k & J^*_{12}\\[0.2cm]J^*_{21} & J^*_{22}-D_2\lambda_k\end{bmatrix};
	\label{ML2}\end{equation}
	\item Linear diffusion with rate $D_P$ defined in \eqref{defDP}:
	\begin{equation}
	\bm{M}_{LP}=\begin{bmatrix} J^*_{11}-D_1\lambda_k & J^*_{12}\\[0.2cm]J^*_{21} & J^*_{22}-D_P\lambda_k\end{bmatrix};
	\label{MLP}\end{equation}
	\item Cross diffusion:   
	\begin{equation}
	\bm{M}_{C}=\begin{bmatrix} J^*_{11}-D_1\lambda_k & J^*_{12}\\[0.2cm]J^*_{21}-J^*_{\Delta 21}\lambda_k & J^*_{22}-J^*_{\Delta 22}\lambda_k\end{bmatrix};
	\label{MC}\end{equation}
	with $J^*_{\Delta 21},\;J^*_{\Delta 22}$ defined in \eqref{JDelta21} and \eqref{JDelta22}.
\end{enumerate}
\medskip 

%\subsection{Turing instability region - linear diffusion}
% Turing INstability Region - linear diffusion
We first want to compare the Turing instability regions of the cases 1 and 2,
% in which we have linear diffusion corresponding to constant diffusion coefficients $D_2,\;D_P$,
 namely the range of the parameters that lead to $\det \bm{M}<0$. The characteristic matrices are \eqref{ML2} and \eqref{MLP} where one should remember that $D_2>D_P$ (because of the modeling assumption). Considering the generic matrix of a (linear) diffusion depending on the parameter $D>0$, we define
$$\bm{M}(D) :=\begin{bmatrix} J^*_{11}-D_1\lambda_k & J^*_{12}\\J^*_{21} & J^*_{22}-D \lambda_k\end{bmatrix} .$$
We compute
%we can study the Turing instability region with respect to $D$ by computing
\begin{align*}
\det \bm{M}(D)&=\left(J^*_{11}-D_1\lambda_k\right)\left(J^*_{22}-D \lambda_k\right)-J^*_{12}J^*_{21}\\
         &=DD_1\lambda_k^2-\left(D_1J^*_{22}+DJ^*_{11}\right)\lambda_k+\det \bm{J}.
\end{align*}
The interesting case (Turing instability appearance) is obtained under the necessary condition $D_1J^*_{22}+DJ^*_{11}>0$, that is 
$$D>\hat{D}:=-\dfrac{D_1J^*_{22}}{J^*_{11}} . $$

%where $\hat{D}$ is a threshold.
%\subsubsection{Solution to \texorpdfstring{$\det J^{*L}_k(D)=0$}{det J^{*L}_k(D)=0}}
The solutions to the equation $\det \bm{M}(D)=0$, which define the boundaries of the Turing instability region, can be written as
\begin{equation}
sol_{1,2} := \dfrac{\left(D_1J^*_{22}+DJ^*_{11}\right)\pm\sqrt{\left(D_1J^*_{22}+DJ^*_{11}\right)^2-4D_1D(J^*_{11}J^*_{22}-J^*_{12}J^*_{21})}}{2D_1D}. \label{eq:sol_12}
\end{equation}
Those solutions exist if the discriminant in \eqref{eq:sol_12} is nonnegative, that is
$$\left(D_1J^*_{22}+DJ^*_{11}\right)^2-4D_1D\det \bm{J}^*\geq 0,$$
which leads to an inequality in $D$:
$$(J^*_{11}D)^2-2D_1(J^*_{11}J^*_{22}-2J^*_{12}J^*_{21})D+(D_1J^*_{22})^2\geq 0.$$
The associated equation has a nonnegative discriminant, so that
$$(D_1(J^*_{11}J^*_{22}-2J^*_{12}J^*_{21}))^2-(D_1J^*_{11}J^*_{22})^2=-4D_1^2J^*_{12}J^*_{21}\det{J_{E_*}}\geq 0 . $$ 
Then the equation admits two real roots
\begin{align*}
\hat{D}_{1,2}&=\dfrac{D_1(J^*_{11}J^*_{22}-2J^*_{12}J^*_{21})\pm 2D_1\sqrt{-J^*_{12}J^*_{21}\det{J_{E_*}}}}{J^{*2}_{11}}\\
             &=\dfrac{D_1}{J^{*2}_{11}}\left[\sqrt{\det{J_{E_*}}}\pm\sqrt{-J^*_{12}J^*_{21}}\right]^2\geq 0,
\end{align*}
which are both nonnegative.

It is also possible to prove that $\hat{D}_1<\hat{D}<\hat{D}_2$. In fact, it can easily be seen that $\hat{D}_2>\hat{D}$ is equivalent to
$$\sqrt{-J^*_{12}J^*_{21}\det \bm{J}^*}>- \det\bm{J}^*,$$
while $\hat{D}_1<\hat{D}$ is equivalent to 
$$\sqrt{-J^*_{12}J^*_{21}\det \bm{J}^*}>\det \bm{J}^* , $$
and can be reduced to $J^*_{21}J^*_{22}<0$. These conditions are therefore always satisfied. Then the values of $D$ which lead to Turing Instability are $D>\hat{D}_2$ (because for $D<\hat{D}_2$, $sol_{1,2}$ are not real or both strictly negative and also for $D<\hat{D}$ no Turing instability can appear).

%\subsubsection{Turing instability region varying \texorpdfstring{$D$}{D}}
We now can perform a qualitative study of the behaviour of the roots $sol_{1,2}(D)$, when we let the parameter $D>\hat{D}_2$ vary. In particular, we have that, from \eqref{eq:sol_12},
$$ \forall D>\hat{D}_2, \qquad sol_{1,2}(D)>0 .$$
Moreover, again from formula \eqref{eq:sol_12}, it can easily be seen that
$$sol_1(\hat{D}_2)=sol_2(\hat{D}_2)=\dfrac{D_1J^*_{22}+\hat{D}_2J^*_{11}}{2D_1\hat{D}_2}>0,$$
and also that
$$\displaystyle{\lim_{D\to +\infty}}{sol_2(D)}=\dfrac{J^*_{11}}{D_1}>sol_2(\hat{D}_2),$$
and 
$$\displaystyle{\lim_{D\to +\infty}}{sol_1(D)}=0.$$
Furthermore, differentiating \eqref{eq:sol_12} with respect to $D$, we obtain
$$\dfrac{\partial}{\partial D}sol_2(D)=-\dfrac{J^*_{22}}{2D^2}+\dfrac{1}{D^2\sqrt{-J^*_{12}J^*_{21}\det \bm{J}^*}}
\left[ \dfrac{J^*_{22}}{2} \dfrac{D_1J^*_{22}+DJ^*_{11}}{D_1D}+\dfrac{\det \bm{J}^*}{D_1D^2}  \right] , $$
and 
$$\dfrac{\partial}{\partial D}sol_1(D)=\dfrac{1}{4D^3D_1}\dfrac{1}{\sqrt{-J^*_{12}J^*_{21}\det \bm{J}}}
\left[ J^*_{22}sol_1(D) -2D\det \bm{J}^* \right] , $$
so that
$$\dfrac{\partial}{\partial D}sol_1(D)>0 , \qquad \dfrac{\partial}{\partial D}sol_2(D)<0.$$ 
This means that the value of $sol_1$ is strictly increasing with respect to $D$, while the value of $sol_2$ is strictly decreasing.
Then, we see that the Turing instability region grows for larger values of $D$. Because of this, the choice of a diffusion rate based only on the behaviour of searching predators would lead to inaccurate conclusions about the possibility of pattern formation.\bigskip

%\subsection{Turing instability region - linear vs cross diffusion}
We then compare the Turing instability regions in the cases 2 and 3, that is when
%  of linear diffusion with constant diffusion rate $D_P$ as defined in \eqref{defDP}, which has
 the characteristic matrices are 
$$\bm{M}_{LP}=\begin{bmatrix} J^*_{11}-D_1\lambda_k & J^*_{12}\\[0.2cm]J^*_{21} & J^*_{22}-D_P\lambda_k\end{bmatrix},$$
and 
%the case of cross diffusion, that is
$$\bm{M}_{LP}=\begin{bmatrix} J^*_{11}-D_1\lambda_k & J^*_{12}\\[0.2cm]J^*_{21}-J^*_{\Delta 21}\lambda_k & J^*_{22}-J^*_{\Delta 22}\lambda_k \end{bmatrix}.$$
We observe that
\begin{align*}
\det \bm{M}_{LP}&=D_1D_P\lambda_k^2-\left(D_PJ^*_{11}+D_1J^*_{22}\right)\lambda_k+\det \bm{J}^*,\\
\det \bm{M}_{C}&=D_1J^*_{\Delta 22}\lambda_k^2-\left(J^*_{\Delta 22}J^*_{11}+D_1J^*_{22}-J^*_{12}J^*_{\Delta 21}\right)\lambda_k+\det \bm{J}^*,
\end{align*}
%and then the Turing Instability regions can be estimated by the study of those determinants.
Both these determinants are second order polynomials in $\lambda_k$ with strictly positive leading coefficients. Furthermore, for $\lambda_k=0$, we know that $\det \bm{M}_{LP}=\det \bm{M}_{C}=\det \bm{J}^*$. We want to compare the leading coefficients $A_L,\;A_C$ of these polynomials on one hand, and the coefficients of $\lambda_k$, that we denote by $B_L,\;B_C$, on the other hand.
 % defined in this way
	%\begin{align*}
	%\det \bm{M}_{LP}&=\underbrace{D_1D_P}_{A_L}k^4-&\underbrace{\left(D_PJ^*_{11}+D_1J^*_{22}\right)}_{B_L}k^2+&\det J^*\\
	%\det J^{*C}_k&=\underbrace{D_1J^*_{\Delta 22}}_{A_C}k^4-&\underbrace{\left(J^*_{\Delta 22}J^*_{11}+D_1J^*_{22}-J^*_{12}J^*_{\Delta 21}\right)}_{B_C}k^2+&\det J^*
%\end{align*}
Those coefficients write:
\begin{align*}
A_L&:=D_1D_P,&\quad& B_L:=D_PJ^*_{11}+D_1J^*_{22},\\
A_C&:=D_1J^*_{\Delta 22}, &\quad&B_C:=J^*_{\Delta 22}J^*_{11}+D_1J^*_{22}-J^*_{12}J^*_{\Delta 21}.
\end{align*}
Substituting the expressions of $D_P$ given in \eqref{defDP} and $J^*_{\Delta 22}$ in \eqref{JDelta22} (in terms of $D_2$ and $D_3$) in $A_L,\;A_C,\;B_L,\;B_C$, we end up with 
the following formulas:
\begin{align*}
A_L&=D_1 \left[ D_2\left( 1-\dfrac{2\mu}{\Gamma}\right) +D_3 \dfrac{2\mu}{\Gamma}\right],\\
A_C&=D_1 \left[ D_2\left( \dfrac{\gamma}{Q_*}+\dfrac{P_*}{Q_*}\left(1-\dfrac{2\mu}{\Gamma}\right)\right) +D_3 \dfrac{2\mu\gamma}{Q_*(\Gamma-2\mu)}\right],\\
B_L&=\left[ D_2\left( 1-\dfrac{2\mu}{\Gamma}\right) +D_3 \dfrac{2\mu}{\Gamma}\right]J^*_{11}+D_1J^*_{22},\\
B_C&=\left[ D_2\left( \dfrac{\gamma}{Q_*}+\dfrac{P_*}{Q_*}\left(1-\dfrac{2\mu}{\Gamma}\right)\right) +D_3 \dfrac{2\mu\gamma}{Q_*(\Gamma-2\mu)}\right]J^*_{11}+D_1J^*_{22}+J^*_{12}\dfrac{D_2-D_3}{Q_*}\left( 1-\dfrac{2\mu}{\Gamma}\right)P_*.
\end{align*}
We first note that 
%look at $A_L$ and $A_C$, coefficients of the leading term of the determinants. B
both $A_L$ and $A_C$ are convex combinations of $D_2$ and $D_3$, since
$$\left( 1-\dfrac{2\mu}{\Gamma}\right)+\dfrac{2\mu}{\Gamma}=1,$$
and
\begin{multline*}
\left( \dfrac{\gamma}{Q_*}+\dfrac{P_*}{Q_*}\left(1-\dfrac{2\mu}{\Gamma}\right) \right)+\dfrac{2\mu\gamma}{Q_*(\Gamma-2\mu)}=\dfrac{1}{Q_*}\left[ \gamma +P_*\left(1-\dfrac{2\mu}{\Gamma}\right)+ \dfrac{2\mu\gamma}{(\Gamma-2\mu)}\right]\\
=\dfrac{1}{Q_*}\left[\dfrac{\Gamma-2\mu}{2\mu}N_*+\dfrac{2\mu\gamma}{(\Gamma-2\mu)} \right] = \dfrac{Q_*}{Q_*}=1.
\end{multline*}
We compare the coefficients of $D_2$ and $D_3$ in those convex combinations. For $D_2$:
$$\dfrac{\gamma}{Q_*}+\dfrac{P_*}{Q_*}\left(1-\dfrac{2\mu}{\Gamma}\right) > 1-\dfrac{2\mu}{\Gamma} \quad \Leftrightarrow \quad N_*>\dfrac{2 \mu \gamma}{\Gamma-2\mu},$$
and for $D_3$:
$$\dfrac{2\mu}{\Gamma} > \dfrac{2 \mu \gamma}{Q_*(\Gamma-2\mu)} \quad \Leftrightarrow \quad N_*>\dfrac{2 \mu \gamma}{\Gamma-2\mu},$$
and those inequalities hold thanks to \eqref{cond4}. As a consequence, we are able to prove that $D_P<J_{\Delta 22}^*$. Indeed:
%\begin{align*}
$$D_2\left( 1-\dfrac{2\mu}{\Gamma}\right) +D_3 \dfrac{2\mu}{\Gamma}<D_2\left( \dfrac{\gamma}{Q_*}+\dfrac{P_*}{Q_*}\left(1-\dfrac{2\mu}{\Gamma}\right)\right) +D_3 \dfrac{2\mu\gamma}{Q_*(\Gamma-2\mu)}.$$
%$$D_2 \left[ \dfrac{\gamma}{Q_*}+\dfrac{P_*}{Q_*}\left(1-\dfrac{2\mu}{\Gamma}\right)-\left( 1-\dfrac{2\mu}{\Gamma}\right)\right]>D_3\left[\dfrac{2\mu}{\Gamma}-\dfrac{2\mu\gamma}{Q_*(\Gamma-2\mu)}\right]$$
%%\end{align*}
%where the coefficients of $D_2$ and $D_3$ are equal (due to the convex combination: $A+B=1$ and $C+D=1$ implies that $C-A=B-D$). The inequality holds because $D_2>D_3$.
We then prove that $B_L>B_C$. In fact, we can write 
\begin{align*}
B_L&=\left[ D_2- (D_2-D_3) \dfrac{2\mu}{\Gamma}\right]J^*_{11}+D_1J^*_{22},\\
B_C&=\left[ D_2- (D_2-D_3) \dfrac{2\mu\gamma}{Q_*(\Gamma-2\mu)}\right]J^*_{11}+D_1J^*_{22}+J^*_{12}\dfrac{D_2-D_3}{Q_*}\left( 1-\dfrac{2\mu}{\Gamma}\right)P_*.
\end{align*}
Starting from these expressions, we have that $B_L>B_C$ if and only if
\begin{multline*}
\left[ \cancel{D_2}- (D_2-D_3) \dfrac{2\mu}{\Gamma}\right]J^*_{11}+\cancel{D_1J^*_{22}}>\\
\left[ \cancel{D_2}- (D_2-D_3) \dfrac{2\mu\gamma}{Q_*(\Gamma-2\mu)}\right]J^*_{11}+\cancel{D_1J^*_{22}}+J^*_{12}\dfrac{D_2-D_3}{Q_*}\left( 1-\dfrac{2\mu}{\Gamma}\right)P_*.
\end{multline*}
Then we can divide by the common strictly positive factor $D_2-D_3$ and multiply both sides by $Q_*$. We obtain
\begin{multline*}
\left[ -\cancel{(D_2-D_3)} \dfrac{2\mu}{\Gamma}\right]J^*_{11}Q_*>
\left[ -\cancel{(D_2-D_3)} \dfrac{2\mu\gamma}{\cancel{Q_*}(\Gamma-2\mu)}\right]J^*_{11}\cancel{Q_*}+J^*_{12}\dfrac{\cancel{(D_2-D_3)}}{\cancel{Q_*}}\left( 1-\dfrac{2\mu}{\Gamma}\right)P_*\cancel{Q_*},
\end{multline*}
and substituting the expressions of $J^*_{11}$ and $J^*_{12}$, we end up with
\begin{multline*}
-\dfrac{2\mu}{\Gamma} \left[r\left(1-\dfrac{2}{\nu}N_* \right)Q_*-\gamma \dfrac{\Gamma-2\mu}{2\Gamma}P_*\right]>
-\dfrac{2\mu\gamma}{(\Gamma-2\mu)Q_*}\left[r\left(1-\dfrac{2}{\nu}N_* \right)Q_*-\gamma \dfrac{\Gamma-2\mu}{2\Gamma}P_*\right]-\dfrac{\gamma^2\mu}{Q_*\Gamma}P_*.
\end{multline*}
We can then expand the product in the r.h.s., and get
$$-\dfrac{2\mu}{\Gamma}\left[r\left(1-\dfrac{2}{\nu}N_* \right)Q_*-\gamma \dfrac{\Gamma-2\mu}{2\Gamma}P_*\right]>-\dfrac{2\mu\gamma}{(\Gamma-2\mu)}r\left(1-\dfrac{2}{\nu}N_* \right).$$
Dividing both sides by $2\mu$, bringing all terms in the l.h.s, we obtain:
$$\left[\dfrac{\gamma}{\Gamma-2\mu}-\dfrac{Q_*}{\Gamma}\right]\left(r-\dfrac{2r}{\nu}N_*\right)+\dfrac{\gamma}{\Gamma}\dfrac{\Gamma-2\mu}{2\Gamma}P_*>0.$$
Using formula \eqref{Qstar} giving $Q_*$ in terms of $N_*$, and eliminating the common factor $\Gamma$, we get
$$\left[\gamma-\dfrac{\Gamma-2\mu}{2\mu}N_*\right]\left(r-\dfrac{2r}{\nu}N_*\right)+\gamma\dfrac{\Gamma-2\mu}{2\Gamma}P_*>0.$$
Using the expression of $P_*$ in terms of $N_*$ in formula \eqref{P*N*}, we obtain
$$\left[\gamma-\dfrac{\Gamma-2\mu}{2\mu}N_*\right]\left(r-\dfrac{2r}{\nu}N_*\right)+\gamma\dfrac{\cancel{\Gamma-2\mu}}{2\cancel{\Gamma}}\dfrac{\cancel{\Gamma}}{2\mu}\dfrac{(\Gamma-2\mu)N_*-2\gamma\mu}{\cancel{(\Gamma-2\mu)}}>0.$$
Using the expression of $N_*$ in formula \eqref{N*pm} (only in the second $N_*$ in the equation above), we end up with the inequality
$$\left[\gamma-\dfrac{\Gamma-2\mu}{2\mu}N_*\right]\left(\cancel{r}-\cancel{r}+\dfrac{\gamma}{2}-\sqrt{\Delta_N}\right)
+\dfrac{\gamma}{4\mu}\bigg((\Gamma-2\mu)N_*-2\gamma\mu\bigg)>0,$$
which is equivalent to 
$$\cancel{\dfrac{\gamma^2}{2}}-\cancel{\dfrac{\gamma}{2}\dfrac{\Gamma-2\mu}{2\mu}N_*}-\sqrt{\Delta_N}\left[\gamma-\dfrac{\Gamma-2\mu}{2\mu}N_*\right]+\cancel{\dfrac{\gamma}{4\mu}(\Gamma-2\mu)N_*}-\cancel{\dfrac{\gamma^2}{2}}>0,$$
and can be reduced to
$$-\sqrt{\Delta_N}\underbrace{\left[\gamma-\dfrac{\Gamma-2\mu}{2\mu}N_*\right]}_{-}>0,$$
which is always true (remember that $P_*>0$).\medskip

\begin{figure}%
\begin{center}
\subfigure[\label{case1}]{\includegraphics[width=.32\columnwidth]{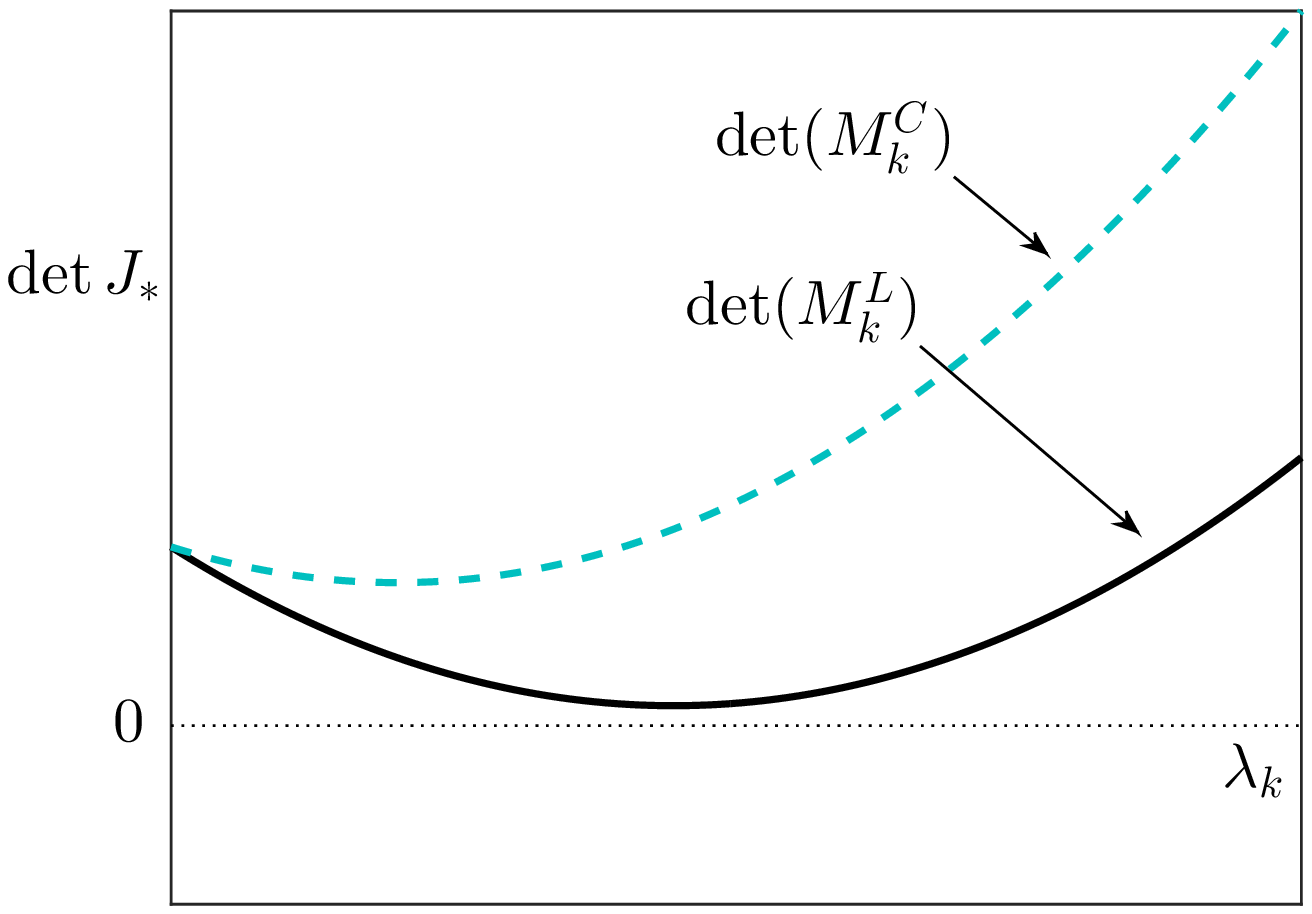}}
\subfigure[\label{case2}]{\includegraphics[width=.32\columnwidth]{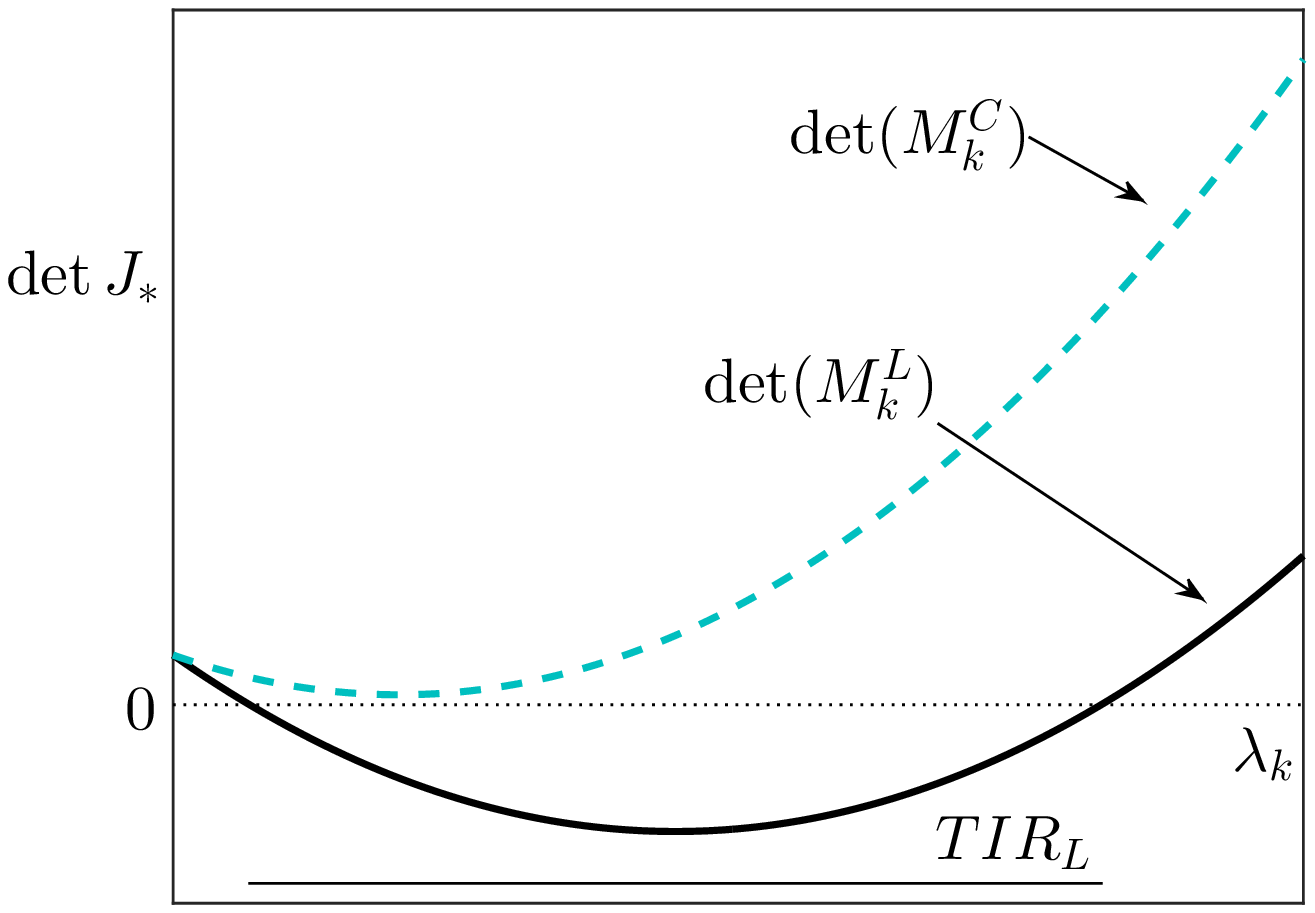}}
\subfigure[\label{case3}]{\includegraphics[width=.32\columnwidth]{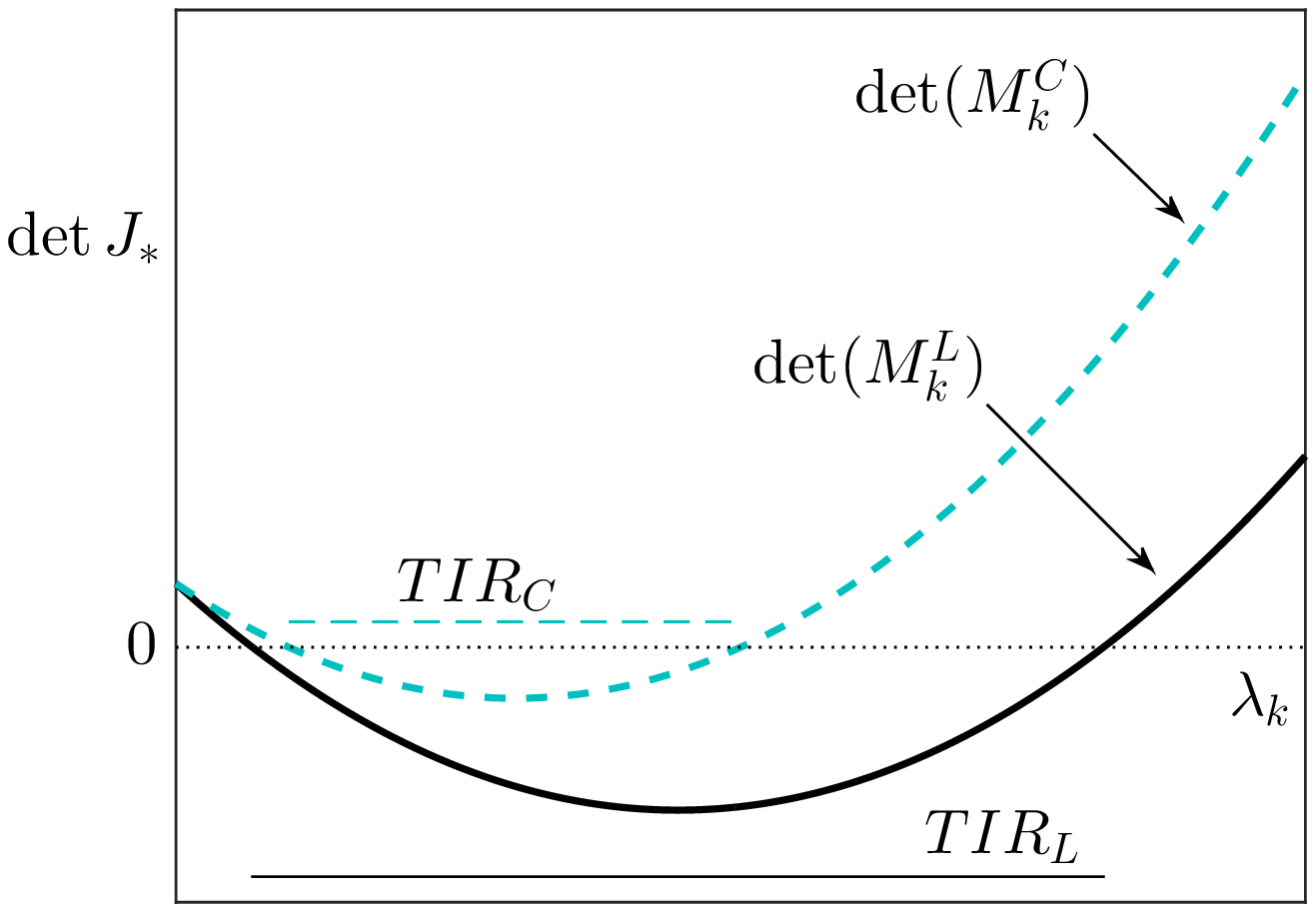}}
\caption{Turing Instability regions for linear diffusion and cross diffusion cases. (a) There are no regions of strictly negative determinant for both linear and cross diffusion, so that in both cases Turing instability cannot appear. (b) The linear diffusion case has a Turing instability region ($TIR_L$), but the determinant of the cross diffusion case is positive for all $\lambda_k$, so that the cross diffusion case does not lead to Turing instability. (c) Both cases lead to nonempty Turing instability regions, but the Turing instability region for the cross diffusion ($TIR_C$) case is strictly included in the Turing instability region of the linear diffusion case ($TIR_L$).}
\label{casesTIR}
\end{center}
\end{figure}

Finally, we see that the determinants of the characteristic matrices with linear and cross diffusion, respectively
$$\det \bm{M}_{LP}=A_L\lambda_k^2-B_L\lambda_k+\det \bm{J}^*, \quad\textnormal{and}\quad	\det \bm{M}_{C}=A_C\lambda_k^2-B_C\lambda_k+\det \bm{J}^* ,$$
and such that $A_C>A_L$ and $B_L>B_C$. Looking at the Turing Instability regions, i.-e. regions in which the determinant of the characteristic matrix is strictly negative, we see that three cases naturally appear:
\begin{enumerate}
	\item There are no regions of strictly negative determinant for both linear and cross diffusion (Figure \ref{case1}).
	\item The linear diffusion case has a Turing instability region, but the determinant of the cross diffusion case is positive for all $\lambda_k$ (Figure \ref{case2}), so that the cross diffusion case does not lead to Turing instability.
	\item Both cases lead to nonempty Turing instability regions (Figure \ref{case3}) and we check that 
	      $$\dfrac{\sqrt{B_L^2-4A_1\det J_*}}{2A_L}>\dfrac{\sqrt{B_C^2-4A_C\det J_*}}{2A_C},$$
				which means that the Turing instability region for the cross diffusion case is strictly included in the Turing instability region of the linear diffusion case.
\end{enumerate}
\bigskip

In all cases, we see that the use of the cross-diffusion model leads to a possibility of obtaining nontrivial patterns which is less likely than when the linear diffusion model is considered. Therefore, the use of a model in which standard diffusion terms are directly added to the reaction terms may lead to an overestimate of the set of the parameters for which patterns appear.

%\clearpage
%\section{Things to do}
%\begin{itemize}
%\item controlla la notazione nella prima parte!!
%\item bisogna spiegre meglio lo scaling!!
%\item comparison between the new BdA and the classical one (for ODEs) in terms of model equilibrium states and shape of the function
% \item significato modellistico aalla cross diffusion che abbiamo derivato?
%\item write down the turing instability analysis for Holling II e cross diffusion.
%\item case $D_3=0$. it shouldn't be a problem for the convergence results, there is something interesting for the Turing instability? 
%\item case $D_2=D_3$?
	%\item Sign of $J^*_{11}$: I proved that there is a threshold under which $J^*_{11}$ is positive. I have still non conditions for the sign of the trace of $ J^*$.
	%\item Comparison with Alonso \cite{Alonso2002}, Lombardo \cite{Tulumello2014}, Huisman \cite{Huisman1997} and other papers \cite{Huang2005,Wang2011,Haque2012}?
	%\item $\Gamma<\gamma$?
	%\item sets of parameters.
	%\item Numerical simulations.
	%\item more involved microscopic model (4 equations)
	%\item forse i pattern ci sono per D3>D2?
%\end{itemize}

\section{Concluding remarks}

This paper focuses on the study of two ``microscopic'' (in terms of time scales) predator-prey models with diffusion, that enable to recover, in a suitable limit, two classical functional responses in the reaction part of the equations and contain a cross-diffusion term. We have also presented rigorous results of convergence of the solutions of these systems towards the solution of the limiting reaction-cross diffusion system.

We first start with two trophic levels, prey and predators, which are further divided into \emph{searching predators} and \emph{handling predators}. The former are predators active in the predation process, the latter are resting individuals. Then, we start from a system of three partial differential equations, with  standard diffusion terms (a constant times the Laplacian), and with a Lotka-Volterra reaction term. Through a \emph{quasi steady-state approximation}, we end up with a system of two PDEs with prey and total predator densities as unknowns, in which a Holling-type II functional response appears together with a cross-diffusion term in the predator equation. This means that the diffusion term relative to predators is much more complicated than a constant times Laplacian of $P$ (linear diffusive term), which in some other models is simply added to the reaction part \cite{Durret}. In particular, the diffusion term obtained in this way  depends on the prey biomass and on both the diffusion coefficients of searching and handling predators $d_2$ and $d_3$. Looking at its expression, the cross-diffusion term reduces the predator diffusion when the prey density increases.

Then we modify the starting model by inserting a competition among predators. With this change we end up after a \emph{quasi steady-state approximation} with a system of two PDEs for prey and total predator densities, characterized by a Beddington-DeAngelis-like functional response, and a cross-diffusion term in the predator equation.
%
%{\tobechecked
%We note that the limiting system presents a functional response close to the Beddington-DeAngelis in the reaction part: in this sense, we have derived this type of functional response by a time-scale argument from a more complicated model. This type of functional response was derived also by Huisman et al. \cite{Huisman1997} in the context of ordinary differential equations, starting from a system of four ODEs by a quasi-steady-state approximation. With their approach, also the logistic growth and the predator turnover terms still depend on this complex expression. Thereafter, they have simplified this complicated quadratic expression with a Pad\'e approximation and they recover the standard formula of the Beddington-DeAngelis functional response.}
%

Also in this case, the limiting system presents a cross-diffusion in the predator equation, which depends on both the diffusion coefficients of searching and handling predators $d_2$ and $d_3$.
%, while the prey diffusion is still linear.  

The Turing instability analysis of the limiting equations is studied in Chapter \ref{SectionTIa}. For the first one, it is known that predator-prey models with a prey-dependent trophic function in the reaction part and (standard) linear-diffusion cannot give rise to Turing instability \cite{Alonso2002}. Even with the cross-diffusion model, no patterns seem to appear under a (biologically reasonable) assumption on the diffusion coefficients. For the second system, in which a Beddington-DeAngelis-like functional response appears, we look for conditions on the parameter values which lead to Turing instability and we compare these Turing instability regions with the ones obtained when the cross-diffusion term is substituted by a standard diffusion. The main point is the fact that the Turing instability region associated to the cross-diffusion system is always strictly included in the Turing instability region of the linear-diffusion system. As a consequence, the use of reaction-diffusion systems for predator-prey interactions of Beddington-DeAngelis type in which standard diffusion is simply added to the reaction terms may lead to an overestimate of the possibility of appearance of patterns (at least in the case when the Beddington-DeAngelis functional response is a consequence of the interactions between searching and handling predators).

It is worth mentioning that in many instances, the introduction of cross-diffusion terms instead of (standard) linear-diffusion terms leads exactly to the opposite result, that is, the increase of the set of parameter values in which patterns develop \cite{Tulumello2014, iida,GambinoLombardo2012}. 
Our study leads then to a rather interesting conclusion: pattern formation originating from Turing instability is counteracted by the cross-diffusion term derived by the Quasi-Steady State Approximation.

\bigskip

L.D. acknowledges support from the French ``ANR blanche'' project Kibord: ANR-13-BS01-0004, and by Universit\'e Sorbonne Paris Cit\'e, in the framework of the ``Investissements d'Avenir'', convention ANR-11-IDEX-0005. C.S. has been partially supported by the French-Italian program Galileo, project G14-34. Support by INdAM-GNFM is also gratefully acknowledged by F.C. and C.S. The authors are very grateful to Odo Diekmann for his valuable suggestions and remarks that improved the manuscript.

\bibliographystyle{elsarticle-num}
\bibliography{bibliography}

%\begin{thebibliography}{10}
%\bibitem{BDF} M. Breden,  L. Desvillettes and K. Fellner.
%
%\bibitem{CDF}
%J. Canizo, L. Desvillettes and K. Fellner.
%\newblock Improved duality estimates and applications to reaction-diffusion equations.
%\newblock {\em Comm. Partial Differential Equations}, 39:1185--1204, 2013.
%
%\bibitem{DFPV} L. Desvillettes, K. Fellner, M. Pierre et J. Vovelle
%
%\bibitem{D_milan} L. Desvillettes
%
%\bibitem{PS} M. Pierre, D. Schmitt
%\end{thebibliography}
\appendix
\end{document}